\documentclass[12pt, twoside]{article}
\usepackage{latexsym,amsfonts,amsthm,amsmath,amscd,amssymb,color,mathrsfs}
\usepackage{graphicx}
\usepackage[linktoc=page,%colorlinks,
linkcolor=blue1,citecolor=green,urlcolor=red,runcolor=red1,linkbordercolor={0 0 1},urlbordercolor={1 0 0},citebordercolor={0 1 0},runbordercolor={1 1 1}]{hyperref}
\usepackage[utf8]{inputenc}
\usepackage{comment}

\newcommand{\con}[2]{\Upsilon^{{#1}\!}{}_{{#2}}}

\newcommand{\Gam}[2]{\Gamma^{{#1}\!}{}_{{#2}}}
\newcommand{\Gamm}[3]{\Gamma_{{#1}\!}{}^{{#2}\!}{}_{{#3}}}

\newcommand{\VV}[2]{V^{{#1}\!}{}_{{#2}}}
\newcommand{\VVV}[3]{V^{{#1}{#2}\!}{}_{{#3}}}

\newcommand{\KKK}[3]{K_{{#1}\!}{}^{{#2}\!}{}_{{#3}}}

\DeclareFontFamily{U}{mathx}{\hyphenchar\font45}
\DeclareFontShape{U}{mathx}{m}{n}{
      <5> <6> <7> <8> <9> <10>
      <10.95> <12> <14.4> <17.28> <20.74> <24.88>
      mathx10
      }{}
\DeclareSymbolFont{mathx}{U}{mathx}{m}{n}
\DeclareFontSubstitution{U}{mathx}{m}{n}
\DeclareMathAccent{\widecheck}{0}{mathx}{"71}

\makeatletter
\let\save@mathaccent\mathaccent
\newcommand*\if@single[3]{%
  \setbox0\hbox{${\mathaccent"0362{#1}}^H$}%
  \setbox2\hbox{${\mathaccent"0362{\kern0pt#1}}^H$}%
  \ifdim\ht0=\ht2 #3\else #2\fi
  }
\newcommand*\rel@kern[1]{\kern#1\dimexpr\macc@kerna}
\newcommand*\widebar[1]{\@ifnextchar^{{\wide@bar{#1}{0}}}{\wide@bar{#1}{1}}}
\newcommand*\wide@bar[2]{\if@single{#1}{\wide@bar@{#1}{#2}{1}}{\wide@bar@{#1}{#2}{2}}}
\newcommand*\wide@bar@[3]{%
  \begingroup
  \def\mathaccent##1##2{%
    \let\mathaccent\save@mathaccent
    \if#32 \let\macc@nucleus\first@char \fi
    \setbox\z@\hbox{$\macc@style{\macc@nucleus}_{}$}%
    \setbox\tw@\hbox{$\macc@style{\macc@nucleus}{}_{}$}%
    \dimen@\wd\tw@
    \advance\dimen@-\wd\z@
    \divide\dimen@ 3
    \@tempdima\wd\tw@
    \advance\@tempdima-\scriptspace
    \divide\@tempdima 10
    \advance\dimen@-\@tempdima
    \ifdim\dimen@>\z@ \dimen@0pt\fi
    \rel@kern{0.6}\kern-\dimen@
    \if#31
      \overline{\rel@kern{-0.6}\kern\dimen@\macc@nucleus\rel@kern{0.4}\kern\dimen@}%
      \advance\dimen@0.4\dimexpr\macc@kerna
      \let\final@kern#2%
      \ifdim\dimen@<\z@ \let\final@kern1\fi
      \if\final@kern1 \kern-\dimen@\fi
    \else
      \overline{\rel@kern{-0.6}\kern\dimen@#1}%
    \fi
  }%
  \macc@depth\@ne
  \let\math@bgroup\@empty \let\math@egroup\macc@set@skewchar
  \mathsurround\z@ \frozen@everymath{\mathgroup\macc@group\relax}%
  \macc@set@skewchar\relax
  \let\mathaccentV\macc@nested@a
  \if#31
    \macc@nested@a\relax111{#1}%
  \else
    \def\gobble@till@marker##1\endmarker{}%
    \futurelet\first@char\gobble@till@marker#1\endmarker
    \ifcat\noexpand\first@char A\else
      \def\first@char{}%
    \fi
    \macc@nested@a\relax111{\first@char}%
  \fi
  \endgroup
}

\newtheorem{lemma}{Lemma}[section]

\newtheorem{rem}[lemma]{Remark}
\newtheorem{prop}[lemma]{Proposition}

\newtheorem{example}[lemma]{Example}
\newtheorem{defn}[lemma]{Definition}

\def\Z{{\mathbb Z}}

\def\L{{\mathcal L}}

\def\E{{\mathbb E}}

\def\E{{\cal E}}

\def\cur{R}
\def\taum{\tau_\M}

\def\h{{\mathfrak h}}

\newcommand{\bea}{\begin{eqnarray}}
\newcommand{\eea}{\end{eqnarray}}

\def\beq{\begin{equation}}
\def\eeq{\end{equation}}

\def\E{{\cal E}}
\def\M{{\cal M}}
\def\L{{\cal L}}

\def\id{{\rm id}}
\def\Vect{{\rm Vect}}
\def\CbracketL{{[[}}
\def\CbracketR{{]]}}

\def\V{{\cal V}}
\def\Tan{{\mathbb T}}

\DeclareMathAlphabet{\mathpzc}{OT1}{pzc}{m}{it}

\def\sk{{\\[-.7em]}}

\begin{document}

\thispagestyle{empty}
\vspace*{-1.5cm}

\begin{flushright}
  {}
\end{flushright}

\vspace{1.5cm}

%%%%%%%%%%%%%%%%%%%%%%%%%%%%%%%%%%%%%%%%%%%%%%%
%%%%%%%%%%%%%%%%%%%%%%%%%%%%%%%%%%%%%%%%%%%%%%%

\begin{center}
{\LARGE On Curvature and Torsion in Courant Algebroids}

\end{center}

%%%%%%%%%%%%%%%%%%%%%%%%%%%%%%%%%%%%%%%%%%%%%%%
%%%%%%%%%%%%%%%%%%%%%%%%%%%%%%%%%%%%%%%%%%%%%%%

\vspace{0.4cm}

\begin{center}
{\bf   Paolo Aschieri$^{1,2,3}$, Francesco Bonechi$^{4}${}\footnote{Corresponding author}, Andreas Deser$^{5}$  }
 
\end{center}

%%%%%%%%%%%%%%%%%%%%%%%%%%%%%%%%%%%%%%%%%%%%%%%
%%%%%%%%%%%%%%%%%%%%%%%%%%%%%%%%%%%%%%%%%%%%%%%

\vspace{0.4cm}

\begin{center}
 $^1$\small{Dipartimento di Scienze  e Innovazione Tecnologica \\
Universit\`a del Piemonte Orientale \\
Viale T. Michel 11, 15121, Alessandria }

  \vspace{0.2cm}

  $^2$\small{Istituto Nazionale di Fisica Nucleare, Sezione di Torino} \\ 
  { Via Pietro Giuria 1, 10125 Torino }\\

  \vspace{0.2cm}

$^3$\small{Arnold--Regge centre, Torino, via P. Giuria 1, 10125, Torino, Italy}

  \vspace{0.2cm}
  
  $^4$\small{Istituto Nazionale di Fisica Nucleare, Sezione di Firenze \\
  Via Giovanni Sansone 1, 50019 Sesto Fiorentino FI}

  \vspace{0.2cm}
  
  $^5$\small{Faculty for Mathematics and Physics, Charles University \\
  Sokolovsk{\'a} 83, 186 75 Praha 8}

  \vspace{0.4cm}

\small{Email: {\tt paolo.aschieri@uniupo.it}, {\tt francesco.bonechi@fi.infn.it}, {\tt andreas3deser@gmail.com} }

\vspace{0.25cm}

\end{center} 

\vspace{0.4cm}

%%%%%%%%%%%%%%%%%%%%%%%%%%%%%%%%%%%%%%%%%%%%%%%
%%%%%%%%%%%%%%%%%%%%%%%%%%%%%%%%%%%%%%%%%%%%%%%

\begin{abstract}
We study the graded geometric point of view of curvature and torsion of $Q$-manifolds (differential graded manifolds).
%Using the definition of a connection in a graded vector bundle over a differential graded manifold, we are able to define the notion of curvature of a graded vector bundle as the square of a covariant derivative as well as the notion of torsion as the derivative of an appropriate identity map. Specifying to the case of degree two differential graded manifolds, 
In particular, we get a natural graded geometric definition of Courant algebroid curvature and torsion, which correctly restrict to Dirac structures. Depending on an auxiliary affine connection $K$, we introduce the $K$-curvature and $K$-torsion of a Courant algebroid connection. These are conventional tensors on the body. Finally, we compute their Ricci and scalar curvature. 
\end{abstract}

%\thispagestyle{empty}

%\eject

\newpage
\tableofcontents

\section{Introduction}
Introduced as the common ambient where pre-symplectic and Poisson structures live \cite{LiuWeiXu}, Courant algebroids are by now a well-studied structure in mathematics and 
physics. Hitchin and Gualtieri's generalized geometry is based on the use of Courant algebroids as a generalization of the tangent bundle. Several basic tools of  differential 
geometry are naturally adapted to this perspective, 
namely the notions of generalized metrics and connections are well established. On the other hand, the definitions of curvature 
and torsion are less straightforward. Indeed the naive generalization obtained by mimicking the standard differential geometric definitions fails to produce 
objects with the expected tensorial properties, due to the failure 
of the Jacobi identity for the Courant bracket.

A proposal for torsion in Courant algebroids was given by 
Gualtieri \cite{Gua2010} and there are several ways to define a Riemann curvature tensor, e.g. by adding additional terms to the naive definition enforcing tensor 
properties \cite{Hohm:2012mf, Jurco:2016emw} or restricting to appropriate subbundles of the generalized tangent bundle \cite{Garcia-Fernandez:2013gja}. The geometric 
interpretation of the Riemann tensor however remains unsatisfactory compared to the ordinary Riemannian curvature. Independently, these ideas turned out to allow 
for a duality-invariant formulation of the bosonic sector of type II supergravity \cite{Coimbra:2011nw, Coimbra:2012af} and play a decisive role in understanding 
contemporary models of duality invariant field theories arising from string theory \cite{Hull:2009mi, Hull:2009zb, Hohm:2010jy}. 

An alternative view on Courant algebroids was established by Roytenberg \cite{Roy2002} and {\v S}evera \cite{Severa:2017oew} in the language of 
graded symplectic geometry: Courant algebroids correspond to degree-2 symplectic differential graded manifolds. 

Later, bundles of the form 
$TM \oplus \wedge^p T^*M$ together with the so called Vinogradov bracket were shown to correspond to symplectic differential graded manifolds of the form $T^*[p+1]T[1]M$ (see \cite{Vin} and \cite{BTZ}). 
Finally, it was realized in \cite{Deser:2014mxa, Deser:2016qkw} that algebraic structures in duality invariant field theories 
inspired by string theory can be understood using symplectic differential graded manifolds of the form $T^*[p]T[1]M$.

The natural definition of connection and curvature we are going to discuss in this paper relies on the following simple fact of algebraic nature. Let 
$({\cal A},Q)$ be a differential graded algebra (dGA) and let $V$ be a graded ${\cal A}$-module endowed with a degree $1$ linear morphism $Q_V:V\rightarrow V$ satisfying
$$
Q_V(a\cdot v) = Q(a)\cdot v + (-)^{\deg a}a \cdot Q_V(v)\;
$$
for each $a\in{\cal A}$ and $v\in V$. The fact that $Q^2=0$ implies that $R_{Q_V}\equiv Q^2_V$ is a degree $2$ ${\cal A}$-module morphism. If $({\cal A}=C(\M),Q)$ is the dGA of
global functions on a graded manifold $\M$ endowed with an homological degree $1$ vector field $Q$ and $V$ is the space of global sections of
a vector bundle $\E\rightarrow\M$, then $Q_V$ and $R_{Q_V}$ are defined by gluing local objects and are the natural generalizations of connection and 
curvature, respectively.

The language of graded geometry is extremely concise and allows to encode rich geometrical structures in terms of few basic notions so that
certain constructions, that in the language of ordinary differential geometry are quite intricate, become canonical and conceptually very simple. The definition of connection on a graded vector bundle is a striking example of this feature and the purpose of this paper is to work out the graded geometrical 
point of view to connection, curvature and torsion for Courant algebroids. We will show that the definitions of curvature and torsion are non ambiguous 
once they are formulated in this language; we will propose also a way to translate them in the language of ordinary differential geometry. We believe that this point of view 
clarifies that the difficulties encountered in the literature contain essential features of this geometry that should be properly taken into account, whatever approach 
one is following.

We introduce first in Section \ref{beginning} the notion of connection and curvature for differential graded manifolds ($Q$-manifolds in short). 
Given $(\M,Q)$ a $Q$-manifold and $\E\rightarrow \M$ a graded vector bundle; a connection on $\E$ is a fibre preserving vector field $Q_\E$ on $\E^*$ projecting  
onto $Q$. Its curvature is simply $Q_\E^2$. 
When $\M$ is non negatively graded (NQ-manifold) we propose a natural definition of torsion. 
The ordinary case corresponds to the case where $\M= T[1]M$, with $M$ a smooth manifold and $Q=d_{dR}$ the de Rham differential on 
forms. In Section \ref{LieAlgbdConn} we illustrate how these definitions coincide with the standard ones in the case of Lie algebroids. In Section 
\ref{mainGraded} we collect the basic needed facts about Courant algebroids. In Section \ref{GradedCourant}, after reviewing the graded 
geometrical description of Courant algebroids, we work out the computations of curvature and torsion. By constructions both have good restriction properties when
we consider a Dirac structure. The objects that we obtain are global sections of some sheaf but
they are not in a natural way sections of some vector bundle over $M$. By choosing an arbitrary affine connection $K$ on $M$, we can define the {\it $K$-curvature} and the 
{\it $K$-torsion} that are sections of vector bundles. This allows the
comparison with other proposals made in a more standard language (see
Section \ref{ClassicalGeometry}); it 
is important to remark that they in general do not contain the full information encoded in the true curvature and torsion. The last two sections are devoted to the study
of these tensors. In particular, in Section \ref{Ricci} we discuss the Ricci tensor and the scalar curvature.

\noindent {\bf{Notations and Conventions}}\\
In the paper we extensively use the notions of graded geometry; we follow the approach of \cite{Roy2002}, see also
the lectures \cite{CS2011}. Here we just fix the notations that we are using. By graded manifold $\M$ we mean a $\Z$-graded sheaf of 
graded commutative algebras over 
the smooth manifold $M$ which is locally isomorphic to $\M(U)\equiv C^\infty(U)\otimes SV$, where $U\subset M$ is an open, $V$ is a graded vector 
space and $SV$ is the graded symmetric algebra of $V$. This local algebra is freely generated and we call the set of 
coordinates on $U$ and generators of $\M(U)$ the coordinates over $U$.  
A global function on $\M$ is a global section of the sheaf and we denote with ${\cal C}(\M)=\oplus_k {\cal C}^k(\M)$ the graded commutative algebra of global functions; a (global) vector field is a derivation of ${\cal C}(\M)$ and we denote with $\Vect(\M)=\oplus_k \Vect^k(\M)$ the graded Lie algebra of such derivations.
For a vector bundle $\E$ over $\M$ we mean a graded manifold $\E$ equipped with an atlas of coordinates 
$\{y^A,s^\alpha\}$ such that $\{y^A\}$ are coordinates for $\M$ and $s^\alpha$
transform linearly; we call  these latter {\it fibre coordinates}. The
dual vector bundle $\E^*$ is defined by the atlas $\{y^A,s_\alpha\}$ where
$s_\alpha$ transforms with the inverse transposed transition functions.  The space of sections of $\E$ is defined as $\Gamma(\E)\equiv C(\E^*)_{lin}\subset C(\E^*)$, {\it i.e.} the 
space of global functions on $\E^*$ which are linear in
the fibre variables.

%Possibly useful:

%$S^\bullet T[1]^*M=\oplus_{n\geq 0} (T[1]^*M)^{\otimes_s n}$ with, as usual, $(T[1]^*M)^0=\mathbb{R}$, 
%and where
%$T[1]^*M=T^*[-1]M\to M$ is the graded vector bundle dual to $T[1]M\to
%M$.\\

%Include also  $\Lambda V$ and $SV$ for $V$ a vector space (needed for
%cooridnate functions algebra on the fibers).

\section{Curvature and Torsion for NQ-manifolds}
\label{beginning}
In differential geometry there are several equivalent ways to define connections and torsion; among them there is one
option mostly suitable to be generalized to the world of graded
geometry. 

Let $M$ be a smooth manifold and $E\rightarrow M$ be a vector bundle; a connection $\nabla$ is a degree $1$ derivation of the differential graded module 
$\Omega(M;E)\equiv\Gamma(\Lambda T^*M\otimes E)$ over
the differential graded algebra (dGA) $(\Omega(M),\wedge, d)$,  {\it i.e.} $\nabla$ satisfies
$$\nabla(\omega\sigma)=d\omega\sigma+(-)^{|\omega|}\omega \nabla(\sigma)\;,$$ 
for each $\omega\in\Omega(M)$ and $\sigma\in\Omega(M; E)$. The {\it curvature of $\nabla$} is 
$\cur_\nabla=\nabla^2$, a $\Omega(M)$-module morphism of degree $2$.

When $E=TM$, then we can pick up
the identity morphism $\id:TM\rightarrow TM$ and look at it  as a 
section $\tau\in\Gamma(T^*M\otimes TM)\subset \Omega(M;TM)$ (in local coordinates $\tau=dx^\mu\otimes \partial_\mu$). We define the {\it torsion} of $\nabla$ as 
\begin{equation}\label{torsion_tangent_bundle}
T_\nabla=\nabla(\tau)\in\Gamma(\Lambda^2 T^*M\otimes TM)\subset \Omega(M;TM).
\end{equation}
The way graded geometry rephrases the above construction is to
consider the $NQ$-manifold $(\M,d)=(T[1]M,d)$ whose dGA of global
functions is $(C(\M),d)=(\Omega(M),\wedge, d)$.
Let $p:T[1]M\rightarrow M$ be the bundle projection and let us
consider the vector bundle $\E\equiv p^*E\rightarrow T[1] M$, and its
space of sections $\Gamma(\E)\equiv
C(\E^*)_{lin}\subset C(\E^*)$.
By introducing a local trivialization of $E$ we can realize the identification
\begin{equation}\label{identification}
\Omega(M;E)\simeq  C(\E^*)_{lin}\subset C(\E^*)~,
\end{equation} 
of the space of sections $\Omega(M;E)$ with the space of sections of $\E$ {\it i.e.}, the space of global functions on the 
dual vector bundle $\E^*\equiv p^*E^*\to T[1]M$ which are linear in
the fibre variables.
%\footnote{In particular setting $E=M\otimes
%\mathbb{R}\to M$ we recover $\Omega({M})=C(T[1]M)$.}
  A connection $\nabla$ on $E$ as above is then the same as the degree $1$ vector field 
$\nabla\in\Vect^1(\E^*)$, preserving
$C(\E^*)_{lin}$ and projecting to $d$ under $\E\rightarrow T[1] M$ 
 ({\it i.e.} the restriction $\nabla|_{C(\M)}$
  leaves invariant $C(\M)$ and it equals $d$). The curvature is the degree $2$ vector field
$$\cur_\nabla=\nabla^2=1/2 [\nabla,\nabla]\in\Vect^2(\E^*)~,$$ preserving
$C(\E^*)_{lin}$ and vertical, {\it i.e.} projecting to the null vector
field on $\M$.
This is the way we express in graded geometry the tensorial property of
the curvature,  indeed the Leibniz rule
$\cur_\nabla(\omega\sigma)=\cur_\nabla(\omega)\sigma+\omega \cur_\nabla(\sigma)$
with verticality $\cur_\nabla(\omega)=0$, imply $\cur_\nabla(\omega\sigma)=\omega\cur_\nabla(\sigma)$ for each $\omega\in\Omega(M)$ and $\sigma\in\Omega(M; E)$.

It is clear that the above definitions of connection and curvature apply to any graded vector bundle $\E\rightarrow\M$; moreover a connection on $p^*E$ defines also a connection on the shifted bundle 
$p^*E[k]$ for any $k\in\Z$.

 In particular let us consider $E=TM$ and
$\E\equiv p^*E[1]= p^*T[1]M\rightarrow T[1] M$. Under the
identification \eqref{identification} the section
$\tau\in\Gamma(T^*M\otimes TM)\subset\Omega(M;TM)$ becomes the
tautological section
%\footnote{The tautological section maps $(m,v_m)
  %\to \taum(m,v_m)=(m,v_m,v_m)$. Given a local frame $\{s_\mu\}$ of
  %$T[1]M$ and the dual frame  $\{s^\mu\}$ of
  %$T[1]^*M$ we identify $\{s^\mu\}$ with fiber coordinate functions on
  %$T[1]M$; then $\taum$ locally reads $\taum=s^\mu\otimes_{C(T[1]M)}
  %s_\mu=s^\mu
  %s_\mu$, indeed
  %$\taum(m,v_m)=(m,v_m,s^\mu(v_m)s_\mu)=(m,v_m,v_m)$. It is useful
%here to think of a section of a (finitely generated) $A$-module as $f^\alpha s_\alpha$ where
%$s_\alpha$ are generating sections, and $f^\alpha\in A={C(T[1]M)}=\Omega({M})$.} 
$\taum\in C(\E^*)_{lin}$ so that the torsion is 
$$
T_\nabla=\nabla(\taum)\in\Gamma(\E) =C(\E^*)_{lin}~.
$$
The shifting of degree is here inessential and simply assigns degree $0$ to 
$\taum$ but it is the natural choice when we consider the general case.\\

Let us introduce local coordinates $\{x^\mu,\psi^\mu\}$ on $T[1]M$,
respectively of degree $(0,1)$. Let $\E\rightarrow \M$ be a vector
bundle, $\{s^\alpha\}$ fiber coordinates of $\E\rightarrow \M$ and 
$\{s_\alpha\}$ fiber coordinates of  $\E^*\rightarrow \M$.
In local coordinates a connection on $\E$ then reads
\begin{equation}\label{conn_local}
\nabla = \psi^\mu\frac{\partial}{\partial x^\mu} + {\Gamma_{\!\mu}}^{\!\alpha\!}{}_\beta \,\psi^\mu
s_\alpha \frac{\partial}{\partial s_\beta}\;.
\end{equation}
Let now  $\E=p^*T[1]M$; in this case we denote the fibre coordinates of degree $1$ as
$s^\mu$, so that the fibre coordinates on $\E^*$ are $s_\mu$ of degree
$-1$.
%\footnote{The identity morphism $\id^*:T[1]^*M \to T[1]^*M$ induces $\id^*:\Gamma(T[1]^*M) \to
%\Gamma(T[1]^*M)$ that maps $s^\mu \to s^\mu=\psi^\mu$ where in the
%last equality we have used the canonical identification
%$\Gamma(T[1]^*M)= C(T[1]M)_{lin}$.}
The tautological section is the
degree zero section 
\begin{equation}\label{torsion_TM}
\taum=\psi^\rho s_\rho\in \Gamma(\E)=C(\E^*)_{lin}
\end{equation} 
that is indeed the image of 
$\tau=dx^\rho\otimes \partial_{\rho}\in \Omega(M;TE)$ under the identification
\eqref{identification}. The torsion reads
$$
T_\nabla = \nabla(\taum) =
{\Gamma_{\!\mu}}^{\!\sigma\!}{}_\rho\,
\psi^\mu\psi^\rho s_\sigma\;.$$ 

\medskip
We are now ready to discuss the general case. Let $\M$ be a $Q$ manifold, {\it i.e.} a graded manifold with a degree 1 cohomological vector field $Q\in \Vect(\M)_1$ satisfying
$Q^2= 1/2 [Q,Q] =0$. Let $\E\rightarrow\M$ be a vector bundle;
%\footnote{
%A  bundle in the $N$-manifolds context is a
%surjection between $N$-manifolds $\E\stackrel{\pi}{\to} \M$. This is a  vector bundle
%if there exists an atlas of local trivializations $\E_{\pi_0^{-1}(U)} \simeq
%\M_U\times \R^{\{k_i\}}$ with transition functions that are linear in the
%fiber coordinates.  The notion of  function linear in the fiber
%coordinates is therefore independent from the trivialization used and hence we
%have the vector space $C^\infty(\E)_{lin}:=\O_{\E\;lin}(E)$ of functions linear in the fibre
%coordinates, where $\O_{\E\;lin}$ is a sheaf of vector spaces, and of
%modules over $\O_{\M}$. The dual sheaf  $\O^\ast_{\E\;lin}$ defined by
%$\O^\ast_{\E\;lin}(U)=(\O_{\E\;lin}(U))^\ast$ defines the dual vector
%bundle $\E^\ast$, with body $E^\ast$ and sheaf of linear functions
%$\O_{\E^\ast\;lin}:=\O^\ast_{\E\;lin}$. The sheaf of algebras
%$\O_{\E^\ast}$ is generated by $\O_{\E^\ast\;lin}$.
%The module  $C^\infty(\E^\ast)_{lin}:=\O_{\E^\ast\;lin}(M)$ is a module over $C^\infty({\cal M}):=\O_{\M}(M)$ and is by
%definition the  $C^\infty({\cal M})$-module of sections $\Gamma(\E)$ of the
%vector bundle $\E$:
%$
%\Gamma(\E)=
%C^\infty(\E^\ast)_{lin}~.
%$
%}
we denote with $C(\E^*)_{lin}\subset C(\E^*)$ the subspace of global functions on the dual vector bundle $\E^*$ 
that are linear in the fibre variables and we call it the space $\Gamma(\E)$ of sections of $\E$.

\begin{defn}\label{def_connection}
A $Q$-connection on $\E\rightarrow \M$ is a degree $1$ vector field $Q_\E\in \Vect^1(\E^*)$ satisfying the following conditions
\begin{itemize}
 \item[-]  $Q_\E$ preserves $\Gamma(\E)$;
 \item[-] $Q_\E$ projects to $Q$.
\end{itemize}
\end{defn}

\begin{rem}{\rm 
A $Q$-connection on $(\M,Q)$  is not a superconnection for the underlying  supermanifold $\M$ (\cite{Quillen1985}). By generalizing 
the discussion at the beginning of this section,
a superconnection on the vector bundle $\E\rightarrow\M$ is a $Q$-connection on $p^*\E\rightarrow T[1]\M$ for the $Q$-manifold $(T[1]\M,d_{dR})$, 
where $d_{dR}$ is the de Rham vector field and $p:T[1]\M\rightarrow\M$.}
\end{rem}

For simplicity, in the rest of the paper we will drop the suffix $Q$ and just refer to $Q_\E$ in Definition \ref{def_connection} as a connection on $\E$.

\begin{defn}\label{def_curvature}
The curvature of a connection $Q_\E$ on $\E\rightarrow\M$ is the degree $2$ vector field $\cur_{Q_\E}=Q_\E^2 = 1/2 [Q_\E,Q_\E]\in \Vect^2(\E^*)$.
\end{defn}

As a consequence of Definition \ref{def_connection}, $\cur_{Q_\E}$ is a
vertical vector field, or equivalently, it is $C(\M)$-linear. 

If the curvature of the $Q$ connection is zero, then in the literature $(\E,Q_\E)$ is known as a {\it Q-bundle} or a {\it dg-vector bundle}, which is a bundle in the category of $Q$-manifolds (see \cite{Mehta2009,Kotov:2007nr,Gruetzmann:2014ica}).

Let us now consider a sub-bundle $\V\subset\E$ over $\M$, {\it i.e.} a
graded embedding $I:\V\rightarrow\E$ respecting the linear structure; 
the dual map $I^*:\E^*\rightarrow \V^*$ is a projection between the
dual vector bundles.% We then have the injections $I:C(\V^*)_{lin}\to C(\E^*)_{lin}$ and $I:C(\V^*)\to C(\E^*)$ and the projections  $I^*:C(\E)_{lin}\to C(\V)_{lin}$ and $I^*:C(\E)\to C(V)$.

\begin{defn}\label{conn_sub_bundle}
We say that a connection $Q_\E$ on $\E$ restricts to a connection $Q_\E|_{\Gamma(\V)}$ on $\V$ if
the degree 1 vector field  $Q_\E\in \Vect(\E^*)^1$ is projectable along $I^*:\E^*\rightarrow \V^*$ to a
vector field on $\V^*$.
\end{defn}

\begin{rem}{\rm 
It is useful to check that the above restriction property is the obvious one when applied to the case of a $Q$-connection for bundles over $T[1]M$, for $M$ smooth manifold. Let $F\rightarrow E$ be an injection of vector bundles over $M$ and let $\nabla$ be a connection on $E$. Let $I:\V=p^*F\rightarrow\E=p^*E$ denote the injection morphism with respect to pull-back bundles along $p:T[1]M\rightarrow M$ 
and let us suppose that the Q-connection $\nabla$ on $\E$ restricts to $\V$ according to Definition \ref{conn_sub_bundle}. The corresponding map $I^*:\E^*\rightarrow \V^*$ can be expressed as $s^a=I_\alpha^a s^\alpha$, where
$\{s^\alpha\}$ and $\{s^a\}$ are fibre coordinates for $p^*E^*$ and $p^*F^*$ respectively. By using (\ref{conn_local}), we compute
$$
\nabla s^a = \psi^\mu (\partial_\mu I^a_\alpha s^\alpha + \Gamma_{\mu\beta}^\alpha s^\beta)=
\psi^\mu \tilde\Gamma^a_{\mu b} s^b$$
where the second equality is true for some $\tilde\Gamma_{\mu b}^a$ since $\nabla$ is projectable along $I^*$.
It is easy now to realize that this condition states that $\nabla$ restricts to the connection on $F$ defined by the coefficients $\tilde\Gamma_{\mu b}^a$.}
\end{rem}

Let now assume that the $Q$ manifold $\M$ is non negatively graded, {\it i.e.} 
all coordinates have non negative degree. We say that $\M$ is a $NQ$
manifold. 
We denote with $M\subset\M$ the body of $\M$. Let $n$ be the degree of $\M$, {\it i.e.} the maximal degree that coordinates of $\M$ can assume.
It is a standard fact (see \cite{Roy2006}) that there exists the following nested fibration of graded manifolds
\begin{equation}\label{nested_fibration}
 \M=\M_n\rightarrow\M_{n-1}\rightarrow\ldots\M_1\rightarrow\M_0\equiv M,
\end{equation}
where $p_j:\M_j\rightarrow \M_{j-1}$ is defined by forgetting all
coordinates of degree $j$ (injecting the sheaf of $\M_{j-1}$ in that
of $\M_j$). In particular, $p_1: \M_1\rightarrow M$ is a vector
bundle. If we denote by 
$p:\M\rightarrow M$ the composition of all the above maps we can
consider the vector bundle $p^*\M_1\rightarrow \M$, and the dual
one $p^*\M_1^*\rightarrow \M$.
%\footnote{{\color{red}The coordinate free definition is not correct. In coordinates if $a=a^\alpha\psi_\alpha$ then you need that 
%$\iota_a=a^\alpha\frac{\partial}{\partial\psi^\alpha}$, which is not global}. Coordinate free definition, to be seen if to be added.
%The tautological section $\taum\in \Gamma(p^*\M_1)=C(p^*\M_1^*)_{lin}$ is then defined by,
%for all $a\in \Gamma(\M_1)=C(\M^*_1)_{lin}$
%\begin{equation}\label{taut_sec2}
%\iota_a\taum=p^*(a)\,
%\end{equation}
%where  $p^*(a)\in \Gamma(p^*\M_1)=C(p^*\M^*_1)_{lin}$ is the pullback
%of $a$, the contraction $\iota_a$ canonically acts on $C(\M_1)_{lin}$
%(since $a\in C(\M^*_1)_{lin}=(C(\M_1)_{lin})^*$ that is the dual
%sheaf) it is then extended as a graded derivation to
%$C(\M_1)$}\footnote{How is it extended to all $C(\M)$? {\color{red}What should be extended ?}}

Let us introduce coordinates 
$\{x^\mu,\psi^\alpha,\ldots\}$ of degree $\{0,1,\ldots\}$ on $\M$. The
fibre coordinates of $\M_1^*\rightarrow \M_0$, dual to the fibre coordinates
$\{\psi^\alpha\}$ of  $\M_1\rightarrow \M_0$, are denoted as
$\{\psi_\alpha\}$ with degree $-1$, their pull back to fibre coordinates 
of the vector bundle $p^*\M_1\to \M$ are denoted by
$s_\alpha=p^*(\psi_\alpha)$, the tautological section 
$\taum\in \Gamma(p^*\M_1)=C(p^*\M_1^*)_{lin}$ in these local
coordinates reads 
\begin{equation}\label{taut_sec}
  \taum = \psi^\alpha s_\alpha\;.
  \end{equation}
%and immediately seen to be globally defined because independent from the choice of coordinates.

\begin{defn}\label{torsion_NQ}
Let $Q_\E$ be a connection on $p^*\M_1\rightarrow\M$. The torsion of $Q_\E$ is the degree $1$ section
$$
T_{Q_\E} = Q_\E(\taum)\in \Gamma(p^*\M_1)\;,
$$
where $\taum$ is the tautological section defined in (\ref{taut_sec}).
\end{defn}

\begin{rem}{\rm
It is clear that if $\M=T[1]M$, the degree of $\M$ is $1$ so the vector bundle $\M_1$ coincides with $T[1]M$ itself and the definition of torsion given in
\ref{torsion_NQ} is the usual one spelled after equation (\ref{torsion_TM}).}
\end{rem}

\bigskip
\bigskip

\section{Lie Algebroid connections}\label{LieAlgbdConn} 
We compare here the usual notion of algebroid connections,
curvature and torsion with the graded geometrical description given in the previous Section.\\

%\noindent{\bf Lie algebroid connection and torsion~}\\
Let $A$ be a Lie algebroid on $M$ with anchor $\rho:A\rightarrow TM$ and Lie bracket on $\Gamma(A)$ denoted as $[\ ,\ ]$. We introduce
local coordinates $\{x^\mu\}$ on $M$ and a local trivialization $\{e_\alpha\}$ on $A$, so that
$$
\rho(e_\alpha) = \rho^\mu_\alpha \partial_\mu\;,\;\;\; [e_\alpha,e_\beta] = f_{\alpha\beta}^\gamma e_\gamma\;.
$$
Let $E\rightarrow M$ be a vector bundle.
A Lie algebroid connection on $E$ is a linear map 
$\nabla:\Gamma(\Lambda^\bullet A^*\otimes E)\rightarrow
\Gamma(\Lambda^{\bullet+1} A^*\otimes E)$ satisfying the Leibniz rule
\begin{equation}\label{leibniz_rule}
\nabla(f\sigma)= d_A(f)\sigma + (-)^{\deg f} f \nabla(\sigma)
\end{equation}
for every $f\in\Gamma(\Lambda^\bullet A^*)$ and $\sigma\in\Gamma(\Lambda^\bullet A^*\otimes E)$, where $d_A$ is the algebroid differential.     
For every section $a\in\Gamma(A)$ and $\sigma\in\Gamma(\Lambda^\bullet
A^*\otimes E)$, using the contraction or inner derivative %\footnote{{\color{red}I would rather skip this and the following footnotes 3,4,5. They are standard computations.}Defined by $i_a(\omega)(a_1,\ldots
%  a_{k-1})=\omega(a,a_1,\ldots a_{k-1})$ for all
%$\omega\in\Gamma(\Lambda^k A^*)$, 
%$k\in \mathbb{N}$. More general notation:
%  $\iota_v\iota_u(\omega\otimes\theta)=\langle v\otimes
%  u,\omega\otimes\theta\rangle=\langle v,\langle
%  u,\omega\rangle\theta\rangle$.}  
$\iota_a: \Gamma(\Lambda^\bullet
A^*\otimes E)\to \Gamma(\Lambda^{\bullet-1}
A^{*}\otimes E)$, we define the covariant derivative $\nabla_a(\sigma)=
(\iota_a\nabla+\nabla \iota_a)(\sigma)$; in particular
$\nabla_a=\iota_a\nabla$ on $\Gamma(E)$.

The curvature of the connection is defined as $F_\nabla=\nabla^2$. Since, as a consequence
of (\ref{leibniz_rule}), $\nabla^2$ is $\Gamma(\Lambda A^*)$-linear, then
$F_\nabla$ is seen as a section of $\Lambda^2 A^*\otimes{\rm End}E$.
If $F_\nabla=0$ then $(E,\nabla)$ is said to be a {\it representation
  of the algebroid $A$}.

Use of the identity $\iota_a\nabla_b-\nabla_b\iota_a=\iota_{[a,b]}$
(or a standard calculation using local trivializations)
shows that
%\footnote{proof:
% $\iota_{a'}\iota_a\nabla^2
%  b=\iota_{a'}(\nabla_a\nabla-\nabla\nabla_a)b=(\nabla_a\iota_{a'}+\iota_{[a',a]})\nabla
%  b-\iota_{a'}\nabla\nabla_a b=([\nabla_a,\nabla_{a'}]-\nabla_{[a,a']})b$.
%} 
%$F_\nabla\in\Gamma(\Lambda^2 A^*\otimes{\rm End}E)$ satisfying, for each $a,a'\in\Gamma(A)$, 
\begin{equation}
F_\nabla(a,a'):=
\iota_{a'}\iota_{a}(F_\nabla)= [\nabla_a,\nabla_{a'}] - \nabla_{[a,a']}\;,
\end{equation}
where $[\;,\,]$ denotes the Lie algebra bracket on $\Gamma(A)$. Finally let $\nabla$ be a connection on $E=A$ itself; the torsion of $\nabla$ is 
$T_\nabla\in\Gamma(\Lambda^2 A^*\otimes A)$ defined as
$T_\nabla=\nabla(\tau_A)$, where $\tau_A\in \Gamma(A^*\otimes A)$ is
the identity section: for all $a\in A$, $\iota_a\tau_A=a$. A standard calculation here too shows that,
for $a, a'\in\Gamma(A)$,%\footnote{The proof uses again the property  $\iota_a\nabla_b-\nabla_b\iota_a=\iota_{[a,b]}$.}
\begin{equation}\label{algebroid_torsion}
T_\nabla(a,a') :=
\iota_{a'}\iota_{a}(T_\nabla) = \nabla_a(a')-\nabla_{a'}(a) - [a,a']~.
\end{equation}
Using
a local frame $\{e_\alpha\}$ of $A$ and its dual $\{e^\alpha\}$ of $A^*$, the
tautological section reads $\tau_A=e^\alpha\otimes e_\alpha$ and  the
torsion reads 
\begin{equation}\label{algebroid_torsion_graded}
  T_\nabla=\nabla(\tau_A)=%d_Ae^a\otimes e_a-e^c\wedge\nabla e_c=f^a_{cd}e^c\wedge e^d\otimes e_a-e^c\wedge R_c^a\otimese_a=
  (\frac{1}{2}f^\alpha_{\beta\gamma}+\Gamma_{\gamma\!}{}_{\beta\!}{}{}^{\alpha\,} )e^\beta\wedge e^\gamma\otimes e_\alpha.\end{equation}
\sk
where the connection coefficients are defined as 
\begin{equation}
\nabla_\eta e^\gamma = 
\Gamma_{\alpha\!}{}_{\beta\!}{}{}^{\gamma\,} \eta^\alpha e^\beta\;,
\end{equation}
with $\eta\in\Gamma(A)$.

Let us switch now to the 
%\noindent{\bf 
$NQ$-manifold description. Let $(A[1],d_A)$ be the $NQ$-manifold encoding the Lie algebroid structure. Indeed, 
the global functions $C(A[1])=\Gamma(\Lambda A^*)$ are a dGA with the algebroid differential $d_A$
seen as a degree $1$ vector field that squares to $0$. If
$\{x^\mu,\xi^\alpha\}$ are the local coordinates of degree $(0,1)$ for
$A[1]$ the differential reads
$$
d_A = \rho^\mu_\alpha \xi^\alpha \frac{\partial}{\partial x^\mu} + \frac{1}{2} f^\alpha_{\beta\gamma}\xi^\beta\xi^\gamma\frac{\partial}{\partial\xi^\alpha}\;.
$$
Let 
$\E\rightarrow A[1]$ be a vector bundle over $A[1]$. Let $Q_\E$ be a connection on $\E$ according to 
Definition \ref{def_connection} and let $F(Q_\E)= Q^2_\E$ be the curvature of $Q_\E$ according to Definition \ref{def_curvature}. If $Q_\E^2=0$ then
$(\E,Q_\E)$ is said to be a {\it representation of the $NQ$-manifold 
$(A[1],d_A)$} (see \cite{Vai1997} where it is called a module over the Lie algebroid). If $\{s^i\}$ are the fibre coordinates for $\E^*$ then
\begin{equation}
\label{Q_algebroid_explicit}
Q_\E = d_A + \Gamma_{i\!\ }{}^j s^i\frac{\partial\ }{\partial s^j}
\end{equation}
where $\Gamma_{i\!\ }{}^j$ is a function of $\{x^\mu,\xi^\alpha\}$.
Since $Q_\E$ is of degree $1$, then $\deg \Gamma_{i\!\ }{}^j=
1-\deg s^i+\deg s^j$. The sections $\Gamma(\E)$
are by definition the functions $C(\E^*)_{lin}$ linear in the fibre coordinates. 
By construction $Q_\E$ preserves $\Gamma(\E)$ and projects to $d_A$.
%\sk

%\noindent{\bf Relation between the two descriptions~} \\
Let us now relate the two definitions of connection and of
algebroid/$NQ$-manifold representation in the case $\E = p^*E$, where
$p:A\rightarrow M$ is the bundle projection and $E$ is a vector bundle over $M$.  Let $\{\sigma^i\}$ be a
local frame (trivialization) of $E$ and $s^i
=p^*(\sigma^i)$ the corresponding local fiber coordinates
of $\E^*$. Since $\deg \Gamma_{i\!\ }{}^j=1$ we have
$\Gamma_{i\!\ }{}^j=\Gamma_{\alpha i\!\ }{}^j\xi^\alpha\,$ that can be seen as the
connection coefficients of an algebroid connection on $E$.
%given
%by
%\begin{equation}
%\nabla_\eta \sigma^\alpha = 
%\Gamma_{a\!}{}^{\alpha\!}{}^{}_{\beta\,} \eta^a \sigma^\beta\;,
%\end{equation}
%with $\eta\in\Gamma(A)$.
%, $\rho:A\rightarrow TM$ denoting the anchor of $A$.

\begin{rem}{\rm
The notion of representation of $NQ$-manifold $(A[1],d_A)$ given above is an important extension of the notion of algebroid representation. Indeed, when 
the vector bundle $\E$ is not concentrated in degree $0$ it corresponds to the concept of {\it representation up to homotopy}, introduced in \cite{AbCr}. For example,
the case $\E=T^*[2]A[1]$ with $Q_\E$ the cotangent lift of $d_A$ is the so called {\it coadjoint representation}, that does not exist in general as an ordinary representation.}
\end{rem}

We finally describe the graded geometric construction of the torsion (\ref{algebroid_torsion}). Since $A[1]$ is of degree $1$ then $\M_1=A[1]$ itself; let $\E=p^*A[1]$ where 
$p:A[1]\rightarrow M$ and let us consider the tautological 
section $\tau\in\Gamma(p^*A[1])=C(p^*A^*[-1])_{lin}$ given by %, for all 
%$a\in \Gamma(A[1])$,   
%$$
%\iota_a\tau=p^*(a) ~,
%$$
$\tau=\xi^\alpha s_\alpha$.
Let $Q_\E$ be a connection on $\E$; the computation of the torsion according to
Definition \ref{torsion_NQ} is exactly the one that leads to (\ref{algebroid_torsion_graded}).

\section{Connections on Courant algebroids}
\label{mainGraded}
We discuss in this section the case of the exact Courant algebroid and its graded geometric interpretation. We recall that a Courant algebroid 
is given by a vector bundle $E\rightarrow M$ equipped with a nondegenerate bilinear form 
$\langle,\rangle$ on the space of sections $\Gamma(E)$,
a bracket $[[,]]:\Gamma(E)\otimes\Gamma(E)\rightarrow\Gamma(E)$ and a bundle map $\rho:\Gamma(E)\rightarrow TM$ called anchor satisfying certain 
compatibility conditions; we refer to \cite{Roy2002} 
for a detailed description. The exact case corresponds to $E=\Tan M=(T\oplus T^*)M$ and the Dorfman bracket 
is defined on sections of $\Tan M$ as
\begin{equation}\label{courant_bracket}
\CbracketL X+\xi,Y+\nu \CbracketR^H = L_X(Y+\nu) - \iota_Y d\xi + \iota_Y\iota_X H 
\end{equation}
where $H$ is a closed three-form. The de Rham cohomology class $[H] \in H^3(M,\mathbb{R})$ is called the {\v S}evera class of the Courant algebroid. 
The bundle $\Tan M$ is endowed with the natural symmetric pairing $\langle X+\xi,Y+\nu\rangle = \nu(X)+\xi(Y)$; the anchor map $\rho:\Tan M\rightarrow TM$ 
is the projection to $TM$.  The bracket (\ref{courant_bracket}) fails to be a Lie algebra bracket since it is not skew-symmetric; this failure is measured by
$$
[[a,a]] =\frac{1}{2} \rho^* d\langle a,a\rangle
$$
for each $a\in\Gamma(\Tan M)$. 

A Dirac structure is a maximally isotropic subbundle $L\subset\Tan M$ that is involutive with respect to the Courant bracket; 
isotropy of $L$ implies that the inherited bracket is skew-symmetric so that it is a Lie algebroid. We denote with $\rho_{T^*}:L\rightarrow T^*M$ the composition of the 
embedding of $L$ into $\Tan M$ with the projection 
to $T^*M$; thanks to isotropy of $L$, $\pi_L=\rho_{T^*}^*\circ\rho:L\rightarrow L^*$ is antisymmetric, {\it i.e.} $\pi_L\in\Gamma(\Lambda^2 L^*)$. Since $L$ is involutive then
$d_L(\pi_L)=0$.

Let $E\rightarrow M$ be a vector bundle. According to \cite{Gua2010} 
a generalized connection is a first order linear differential operator $D:\Gamma(E)\rightarrow \Gamma(\Tan M\otimes E)$ such that
\begin{equation}\label{gualtieri_connection}
D(f\sigma) = f D\sigma + (\rho^*df)\otimes\sigma\;,
\end{equation}
for $f\in C^\infty(M)$ and $\sigma\in\Gamma(E)$. We denote $D_a = \iota_a D$ for $a\in\Gamma(\Tan M)$. It is easy to check that $D=\nabla + V$ where
$\nabla$ is an ordinary connection on $E$ and $V\in\Gamma(TM\otimes{\rm End} E)$. The curvature operator $R_D$ is defined as 
\begin{equation}\label{gualtieri_curvature_operator}
R_D(a,b) = [D_a,D_{b}] - D_{[[a,b]]_{sk}}
\end{equation}
for each $a,b\in\Gamma(\Tan M)$ and $[[a,b]]_{sk}=\frac{1}{2}([[a,b]]-[[b,a]])$. Since the bracket (\ref{courant_bracket}) is not skewsymmetric, $R_D(a,b)$ is not tensorial in $a,b$,  indeed, $R_D(f a, b) = f\,R_D(a,b) + \langle a, b\rangle df$ for $a,b \in \Gamma(\Tan M)$ and $ f\in C^\infty(M)$. Note that the restriction of $R_D$ to a Dirac structure $L$ is tensorial, since $\langle a, b\rangle$ vanishes if $a,b \in \Gamma(L)$.
In fact if $D$ preserves $L$, the restriction of $D$ to $L$ is an algebroid connection with curvature the restriction of $R_D$.

It is common in the literature (see e.g. \cite{Jurco:2016emw}) to demand compatibility with the bilinear symmetric pairing of a Courant algebroid connection $D$, i.e.
\begin{equation}\label{compatibility_pairing}
  \rho(a)(\langle b,c\rangle) =\,\langle D_a b, c\rangle + \langle b, D_a c\rangle \;, \quad a,b,c \in \Gamma(E)\;.
\end{equation}
We will not assume this condition so far but it will be needed in section \ref{Ricci}.

Let us choose $E=\Tan M$ and let $D$ be a generalized connection on $\Tan M$. In \cite{Gua2010} the torsion ${\mathfrak T}_D\in\Gamma(\Lambda^2 \Tan M\otimes \Tan M)$ of $D$ is defined by
\begin{equation}\label{gualtieri_torsion}
{\mathfrak T}_D(a,b,c) = \langle D_ab-D_b a - [[a,b]]_{sk},c\rangle + \frac{1}{2} (\langle D_ca,b\rangle - \langle D_c b, a\rangle)
\end{equation}
where $a,b,c\in\Gamma(\Tan M)$ . Even when $D$ preserves a Dirac structure $L$, {\it i.e.} 
$D_ab\in\Gamma(L)$ for $a,b\in\Gamma(L)$,
${\mathfrak T}_D$ does not restrict in general to the Lie algebroid torsion for $L$ defined in (\ref{algebroid_torsion}). 
In contrast to the notion of curvature operator 
(\ref{gualtieri_curvature_operator}), that is tensorial only when restricted to a Dirac structure, the torsion (\ref{gualtieri_torsion}) is tensorial on the whole Courant algebroid.
%
%{\color{red}We will be interested in the more general case of an exact Courant algebroid, where the Dorfman bracket is twisted by a closed three form $H$
%\begin{equation}
%  [[X+\xi, Y+\nu]]^H = L_X(Y + \nu) - \iota_Y d\xi + \iota_Y\iota_X H\;.
%\end{equation}
%The bracket $[[\cdot,\cdot]]^H$ is usually termed $H$-twisted Dorfman bracket. 
%The de Rham cohomology $[H] \in H^3(M,\mathbb{R})$ is called the {\v S}evera class of the Courant algebroid}. 
%The form $H$ will appear as twist of the cohomological vector field corresponding to a Courant algebroid in section \ref{GradedCourant}. 

Finally, a \emph{generalized Riemannian metric} on $M$ is an orthogonal, self-adjoint bundle map $\h \in \textrm{End}(\Tan M)$ with the property that 
$\langle \h(a), a\rangle , a\in \Gamma(\Tan M)$ is positive definite \cite{Gualtierithesis}. This means that $\h^2 = \textrm{id}$ and thus it gives 
rise to $\pm 1$ eigenbundles of $\Tan M$. It is well-known (\cite{Gualtierithesis}) that the positive eigenbundle can be expressed as the graph of the sum of 
a metric $g$ and a two-form $B$ on $M$; it is customary to represent it via the symmetric bilinear form  $G(a,b):=\langle \h(a),b\rangle$ on $\Tan M$  given by the 
symmetric matrix
\begin{equation}\label{GeneralizedMetric}
  G = \begin{pmatrix}
    g - B\,g^{-1}B & B\,g^{-1}\\
    -g^{-1}B & g^{-1}
  \end{pmatrix}\;.
\end{equation}
We will use this explicit form of a generalized Riemannian metric in section \ref{Ricci} to define and compute the generalized Ricci curvature.

\section{Graded geometry of Courant Algebroids}
\label{GradedCourant}

From \cite{Roy2002} we know that the structure of a Courant algebroid can be encoded in the data of a 2-symplectic NQ manifold. In the case of an exact Courant algebroid, 
the construction goes as follows. Let us consider the 
NQ manifold $(\M= T^*[2]T[1]M,d_\M)$. Indeed, if $(x^\mu,\psi^\mu,b_\mu,p_\mu)$ are the local coordinates of degree $(0,1,1,2)$, respectively, 
then $d_\M$ is the degree $1$ cohomological vector field defined as
\begin{displaymath}
d_\M\,=\,\psi^\mu\frac{\partial}{\partial
  x^\mu}+\Bigl(p_\mu + \frac{1}{2} H_{\mu \nu\rho} \psi^\nu \psi^\rho\Bigr)\frac{\partial}{\partial b_\mu} - \frac{1}{3!} \partial_\kappa H_{\mu\nu\rho}\psi^\mu\psi^\nu \psi^\rho\frac{\partial}{\partial p_\kappa}~.
%d_\M x^\mu = \psi^\mu\;,\;\;\; d_\M b_\mu = p_\mu\;.
\end{displaymath}
Moreover $\M$ is canonically 2-symplectic, with symplectic form given by $\omega= dx^\mu dp_\mu + d\psi^\mu d b_\mu$ so that $d_\M$ is hamiltonian with hamiltonian
$\Theta\in C^3(\M)$ defined by 
\begin{displaymath}
\Theta = \psi^\mu p_\mu + \frac{1}{3!} H_{\mu \nu\rho} \psi^\mu\psi^\nu\psi^\rho\;,
\end{displaymath}
that satisfies $\{\Theta,\Theta\}=0$ if and only if $dH = 0$. %\footnote{By degree reasons, $\omega$ is exact and there is a canonical primitive. Indeed,
%let $\epsilon= \psi^\mu\partial_{\psi_\mu}+b_\mu \partial_{b_\mu}+2 p_\mu\partial_{p_\mu}$ the Euler vector field then $L_\epsilon \omega = 2 \omega$ implies that $\omega = d\theta$ where
%begin{equation}\label{primitive_omega}
%\theta = \frac{1}{2} \iota_\epsilon \omega = \frac{1}{2}\psi^\mu db_\mu + \frac{1}{2} b_\mu d\psi^\mu - p_\mu %dx^\mu
%\end{equation}}
The transformation of coordinates can be easily checked to be the following ones:
% if $x,\psi,b, p$ are the column vectors of local
%coordinates of $T^*[2]T[1]M$, then 
under a coordinate change in $M$,
$x'^\mu=\varphi^\mu(x)$, we have $\psi'=J\psi$ and $b'={J^{-1}}^Tb$
with matrix $J=(\partial\varphi^\mu/\partial x^\nu)$, while
$p'={J^{-1}}^Tp-{J^{-1}}^TdJ^T{J^{-1}}^Tb$ with $dJ=\partial_\mu J\,\psi^\mu$.

Let $\E\rightarrow \M$ be a graded vector bundle over $\M$ and let us consider a connection $Q_\E$ on $\E$ as in Definition \ref{def_connection}. 
Let $\{s_\alpha\}$ be the local fibre coordinates of $\E^*$; we then have that
\begin{equation}
\label{Q_explicit}
Q_\E = d_\M + \con{\alpha}{\beta} s_\alpha\frac{\partial\ }{\partial s_\beta}\;.
\end{equation}

The curvature  $\cur_{Q_\E}$  of the connection $Q_\E$ is, as in Definition \ref{def_curvature}, the vertical vector field $Q_\E^2= 1/2 [Q_\E,Q_\E]$. An explicit computation shows that
\begin{equation}
\label{F_explicit}
\cur_{Q_\E} = (d_\M\con{\alpha}{\beta} + (-)^{|s_\alpha|+|s_\beta|\,} \con{\alpha}{\gamma} \con{\gamma}{\beta}) s_\alpha\frac{\partial\ }{\partial s_\beta}\; .
\end{equation}

%Since our aim is to compare with the definitions considered in the
%classical approach to the differential geometry on Courant algebroids
We  now discuss  the special case $\E= p^*E[k]$, for $k\in\Z$, where
$p:\M\rightarrow M$ is defined by $p(x^\mu,\psi^\mu,b_\mu,p_\mu)=x^\mu$ and $E\rightarrow M$
is a vector bundle on $M$. In this case the fibre coordinates $s^\alpha$ have degree $k$ and $|\con{\alpha}{\beta}|=1$ so that we can expand
\begin{equation}\label{R_explicit_non_graded_vb}
\con{\alpha}{\beta} =\Gam{\alpha}{\beta} + \VV{\alpha}{\beta} = \Gamm{\mu}{\alpha}{\beta\,} \psi^\mu +\VVV{\mu}{\alpha}{\beta\,}  b_\mu\;.
\end{equation}
By looking at the transformation properties of $\{\con{\alpha}{\beta}\}$ in
(\ref{Q_explicit}), one checks that ${\Gamma_\mu}^\alpha{}_\beta$ are
the connection coefficients of a connection $\nabla^\Gamma$ on $E$ and
${V^\mu}^\alpha{}_\beta$ the components of a tensor
$V\in\Gamma(TM\otimes {\rm End}E)$. We then proved the following proposition.

\begin{prop}A $Q$-connection on $p^*E[k]$ for $(T^*[2]T[1]M, d_\M)$ is equivalent to a generalized
connection $D=\nabla^\Gamma+V$ on $E$ as defined in (\ref{gualtieri_connection}).
\end{prop}

After a straightforward computation
%\footnote{ ????? it seems
%  rather short and straighforward, just recall
%$d_\M=d+p_\mu\frac{\partial}{\partial b_\mu}$ and use the matrix
%notation $\Gamma+V=\Gam{\alpha}{\beta}+\VV{\alpha}{\beta}$, then 
%$$
%\cur_{Q_\E}^{~~\alpha\!}{}_\beta =d\Gamma+\Gamma\Gamma+V V+ dV+\Gamma
%V +V\Gamma + (d_\M-d)V
%=R_{\nabla^\Gamma} +  V V+ \nabla^\Gamma V+   (d_\M-d)V$$}
the curvature reads
\begin{eqnarray}\label{allthecomponents}
(\cur_{Q_\E})^{\alpha\!}{}_\beta &=& (\partial_\mu \Gamm{\nu}{\alpha}{\beta}+
\Gamm{\mu}{\alpha}{\gamma}\Gamm{\nu}{\gamma}{\beta} + \tfrac{1}{2}H_{\rho \mu\nu}\,V^\rho{}^\alpha{}_\beta)\psi^\mu\psi^\nu +(\partial_\mu \VVV{\nu}{\alpha}{\beta}+ \VVV{\nu}{\gamma}{\beta}
                                     \Gamm{\mu}{\alpha}{\gamma} 
\cr
& & -\Gamm{\mu}{\gamma}{\beta} 
\VVV{\nu}{\alpha}{\gamma})\psi^\mu b_\nu-\VVV{\mu}{\alpha}{\gamma}\VVV{\nu}{\gamma}{\beta}\, b_\mu b_\nu + 
\VVV{\mu}{\alpha}{\beta}\,p_\mu ~.
\end{eqnarray}
Therefore we finally obtain the following formula for the curvature
\begin{equation}\label{F_explicit_non_graded_vb}
R_{Q_\E}=R_\nabla +  V\, V + R_{Q_\E}^{(1,1)}\,,
\end{equation}
where $R_\nabla$ is the curvature of the connection $\nabla^\Gamma$ and
\begin{equation}\label{F_explicit_non_graded_vb11}
(R_{Q_\E}^{(1,1)})^{\alpha\!}_{\, \beta} = (\partial_\mu \VVV{\nu}{\alpha}{\beta}+ \VVV{\nu}{\gamma}{\beta}
                                     \Gamm{\mu}{\alpha}{\gamma}-\Gamm{\mu}{\gamma}{\beta} 
\VVV{\nu}{\alpha}{\gamma})\psi^\mu b_\nu + 
\VVV{\mu}{\alpha}{\beta}\,p_\mu ~.
\end{equation}
The  term $R_{Q_\E}^{(1,1)} $  is not a section of a vector bundle on $M$ due to the
transformation properties of the coordinates $p_\mu$.  We remark that $p^*E$ is a $Q$-bundle if and only if $V=0$ and
$\nabla$ is a flat connection.

We now discuss the torsion. The degree of $T^*[2]T[1]M$ is $2$ and the
fibration of graded manifolds described in \eqref{nested_fibration} in
this case reads 
$$
T^*[2]T[1]M\rightarrow \Tan[1]M\rightarrow M\;,
$$
so that $\M_1\rightarrow M$ is $\Tan[1]M\rightarrow M$. Let us consider then a connection $Q_\E$ on $\E\equiv p^*\M_1$; it encodes a generalized connection $D=\nabla + V$ where $\nabla$ is a connection on $\Tan M$ and $V\in\Gamma(T^*M\otimes{\rm End}(\Tan M))$. 
Let us introduce the fibrewise coordinates $\{s^\mu,s_\mu\}$ of degree $-1$ of $\E^*$; the tautological section is then
$$
\taum = \psi^\mu s_\mu + b_\mu s^\mu\in \Gamma(\E)\subset C(\E^*)\;.
$$
As in Definition \ref{torsion_NQ}, the torsion of $Q_\E$ is then 
given by
$$
T_{Q_\E} = Q_{\E}(\taum)\;.
$$
A straightforward computation, see Appendix \ref{appA} for details on the conventions and the result in components, shows that
\begin{equation}
\label{torsion_components}
T_{Q_\E} = T_\nabla+ T_V + T_{Q_\E}^{(1,1)}
\end{equation}
where $ T_{Q_\E}^{(1,1)} \in \Gamma(\E)$, given in (\ref{torsion_mixed_term}), is 
not a section of  a vector bundle on $M$, rather  it is a global
section of the sheaf $\E^*$ on $M$.
On the other hand $T_V\in \Gamma(\Lambda^2TM\otimes \Tan M)$, and is defined,  for all  $\nu,\lambda\in \Gamma(T^*M)$, as 
\begin{equation}
\langle T_{V}, \nu\wedge \lambda\rangle=V_\nu \lambda-V_\lambda \nu~,
\end{equation}
and $T_\nabla\in\Gamma(\Lambda^2T^*M\otimes \Tan M)$ is defined for all $X,Y\in \Gamma(TM)$ as
\begin{eqnarray}
\langle T_\nabla,X\wedge Y\rangle=\nabla_X Y-\nabla_Y X -[X,Y] ~.
\end{eqnarray}
In order to analyze this latter expression we decompose the 
connection $\nabla$ on $\Tan M$, entering $D=\nabla^\Gamma +V$
as follows 
\begin{equation}
\label{nabladecomp}
\nabla^\Gamma = \left( \begin{array}{cc} \nabla^{TT} & \nabla^{TT^*}\cr
\nabla^{T^*T}&\nabla^{T^*T^*}
\end{array}\right)\;,
\end{equation}
where $\nabla^{TT}$ and $\nabla^{T^*T^*}$ are connections on $TM$ and $T^*M$, while $\nabla^{TT^*}\in\Gamma(T^*M\otimes 
{\rm Mor}(T^*M,TM))$ and $\nabla^{T^*T}\in\Gamma(T^*M\otimes {\rm Mor}(TM,T^*M))$.
Then we have 
\begin{equation}
T_\nabla=T_{\nabla^{TT}}+T_{\nabla^{T^*T}}\;,
\end{equation}
where $T_{\nabla^{TT}}$ is the torsion of the connection on $TM$, while
$T_{\nabla^{T^*T}}\in \Gamma(\Lambda^2T^*M\otimes T^*M)$ is defined by
\begin{equation}
\langle T_{\nabla^{T^*T}}, X\wedge Y\rangle=\nabla_X^{T^*T} Y-\nabla_Y^{T^*T} X + H(X,Y,\cdot)~. 
\end{equation}

\subsection{$K$-curvature and $K$-torsion}

The curvature and torsion that we computed above in (\ref{F_explicit_non_graded_vb}) and (\ref{torsion_components}) are global objects but are not sections of some vector bundle over $M$. Namely, $\M$ is not a graded vector bundle over $M$ and so the global functions ${\cal C}(\M)$ are not sections of any vector bundle in a canonical way. This
is a direct consequence of the nature of the geometry encoded in the
Courant algebroid.  Nevertheless, $\M$ can be made a graded vector bundle over $M$ in a non canonical way by introducing an auxiliary affine connection  
$\nabla^K$ on $M$ that provides the following splitting 
\begin{equation}\label{splitting_courant}T^*[2]T[1]M~\simeq~ T[1]M\oplus T^*[1]M\oplus T^*[2]M~.\end{equation}
By using this identification we will be able to express all geometrical objects in a more friendly form at the cost
of introducing an auxiliary object. We are going to discuss here this point of view since it will be useful for later comparison with the other approaches in literature. 

The splitting (\ref{splitting_courant}) can be seen as the following change of variables in $\M$:  let $K_\nu{}^\rho{}_{\!\mu}$ denote the connection coefficients of $\nabla^K$,  we then define
\begin{equation}\label{ptilde}
\tilde p=p+K^Tb~~~i.e.~~~~~\tilde{p}_\mu=p_\mu + K_{\nu}{}^{\rho}{}_{\!\mu} \psi^\nu b_\rho\;.
\end{equation}
It is easy to see that $\tilde{p}_\mu$ transforms as a
tensor on $M$, $\tilde p'={J^{-1}}^T\tilde p_{\,}$; (indeed $K'=JKJ^{-1}-JdJ^{-1}$, with matrix $K=(K^\rho{}_\sigma$), with $K^\rho{}_\sigma=K_{\nu}{}^{\rho}{}_{\!\sigma} \psi^\nu$). We can now express curvature and torsion in the new
variables $(x^\mu,\psi^\mu,b_\mu,\tilde p^\mu)$.

The term $R_{Q_\E}^{(1,1)} $ in (\ref{F_explicit_non_graded_vb}) can now be read as the sum of two tensors, namely
\begin{eqnarray}\label{F_explicit_non_graded_vb11_bis}
(R_{Q_\E}^{(1,1)})^{\alpha\!}_{\, \beta} &=& (\partial_\mu \VVV{\nu}{\alpha}{\beta}+ \VVV{\nu}{\gamma}{\beta}
                                     \Gamm{\mu}{\alpha}{\gamma}-\Gamm{\mu}{\gamma}{\beta} 
\VVV{\nu}{\alpha}{\gamma} - \VVV{\rho}{\alpha}{\beta} K_{\mu}{}^{\nu}{}_{\!\rho})\psi^\mu b_\nu + 
\VVV{\mu}{\alpha}{\beta}\,\tilde{p}_\mu \cr
&=&  \nabla^{\Gamma K}_\mu\VVV{\nu}{\alpha}{\beta}\psi^\mu b_\nu + \VVV{\mu}{\alpha}{\beta}\,\tilde{p}_\mu~.
\end{eqnarray}
The first summand of (\ref{F_explicit_non_graded_vb11_bis}) is 
$\nabla^{\Gamma K}(V) \in 
\Gamma(T^*M\otimes TM \otimes {\rm End}(E))
$, the covariant derivative of $V$ with respect to the connection 
$\nabla^{\Gamma K}$ on the tensor product bundle $TM\otimes {\rm End}(E)$. 

We define the {\it $K$-curvature tensor} of the generalized connection $D=\nabla^\Gamma + V$ on $E$ as the following tensor
\begin{equation}\label{K_curvature}
R^K_D = R_\nabla  +  V\wedge V + \nabla^{\Gamma K}(V) \in\Gamma(\Lambda^2\Tan M\otimes{\rm End}E)
\end{equation}
so that 
\begin{equation}
\label{Full_vs_Kcurv}
R_{Q_{\cal E}} = R^K_D + V^\mu\tilde{p}_\mu.
\end{equation}

%From the above computations we get
%\begin{equation}
%F^K(a,b) = F(\Gamma)(X,Y) - \frac{1}{2}[V,V](\xi,\eta) + (\nabla_X^{\Gamma K}V)(\eta)-(\nabla_Y^{\Gamma K}V)(\xi)\;,
%\end{equation}
%and $a=X+\xi$, $b=Y+\eta$. We call $F^K$ the {\it $K$-curvature tensor of $D$}.\\

We can analogously proceed with the torsion computed in \eqref{torsion_components}. The non tensorial term
in \eqref{torsion_components} is $T^{(1,1)}_{Q_\E}$, that, after the change of variables, can be written as
$$
T^{(1,1)}_{Q_\E} = T^{(1,1)K}_{\nabla,V} +\tilde{p}_\mu s^\mu
$$
where $T^{(1,1)K}_{\nabla,V}\in\Gamma(TM\otimes T^*M\otimes\Tan M)$ is defined as
\begin{equation}\label{tensorial_torsion_3}
T^{(1,1)K}_{\nabla,V} (X,\nu) =\nabla_X(\nu)-\nabla^K_X(\nu)-V_\nu(X) 
 %\cr & &+ \langle \nabla_X^{TT^*}(\nu)- V_\nu^{TT}(X), c_{T^*}\rangle
\;\;,
\end{equation}
where $X\in \Gamma(TM)$ and $\nu\in\Gamma(T^*M)$. We finally define the {\it $K$-torsion tensor} of the generalized connection $D=\nabla^\Gamma + V$ on $\Tan M$ as
\begin{equation}\label{K_torsion}
T_D^K = T_\nabla+ T_V + T^{(1,1)K}_{\nabla,V} \in\Gamma(\Lambda^2\Tan M\otimes \Tan M)
\end{equation}
so that
\begin{equation}\label{Full_vs_Ktor}T_{Q_\E}= T_D^K + \tilde{p}_\mu s^\mu.
\end{equation}

\begin{rem}\label{second_def_K_tensor}{\rm
The $K$-curvature and torsion above can be introduced in the following equivalent way, by using contraction operators. We recall that the contraction operator 
$\iota_v$  along $v\in\Gamma(TM)$ acting on forms is seen as a vector field $\iota_v\in Vect(T[1]M)$; given any vector bundle $E$, $\iota_v$ lifts canonically to a vector field on 
$p^*E^*$, preserving $\Gamma(E\times\Lambda T^*M)$ seen as the linear functions $C(p^*E^*)_{lin}$. 

In the case of a graded vector bundle $\E$ over $\M$, we first observe that
every vector field $v\in Vect(\M)$ that has zero component in the direction
$\partial/\partial x^\mu$ can be lifted to $\E$. We will be interested in vector fields defined from sections of the generalized tangent bundle $\Tan M$.
Indeed for every $a=a^\mu\partial_\mu+a_\mu dx^\mu\in\Gamma(\Tan M)$ we can define the $K$ dependent contraction vector field $\iota_{a}^K\in Vect(\M)$ as
\begin{equation}\label{contraction_vf}
\iota_{a}^K = a^\mu(\frac{\partial\ }{\partial\psi^\mu}+
\KKK{\mu}{\rho}{\nu} b_\rho\frac{\partial\ }{\partial p_\nu})+
a_\mu(\frac{\partial\ }{\partial b_\mu}-
\KKK{\rho}{\mu}{\nu} \psi^\rho\frac{\partial\ }{\partial p_\nu}) \;.
\end{equation}
The correction depending on the connection coefficients $\KKK{\rho}{\mu}{\nu}$ is needed in order to define a global vector field of $\M$, again as a consequence of the transformation properties (\ref{ptilde}) of $p_\mu$.
It is now easy to see that for every $a,b\in\Gamma(\Tan M)$ we have
$$
R^K_D(a,b) =[\iota^K_b,[\iota^K_a, R_{Q_{{\cal E}}}]]\;,\;\;\;\; T_D^K(a,b) = \iota^K_b\iota^K_a T_{Q_\E}\;.
$$ 
It is clear that the rationale of the above formula is that $R_{Q_\E}$ is a vector
field.}
\end{rem}

\subsection{Dirac structures}\label{obstr_dirac}

Let $L\subset \Tan M$ be a Dirac structure. In the graded language, this is equivalently described as a $d_\M$-invariant lagrangian submanifold $\L_L\subset\M$. 
See the Appendix \ref{appB} for an explicit description. 
%such that $\theta|_{\L_L}=0$, where $\theta$ is defined in \ref{primitive_omega}. 
Since, according to Definition \ref{def_connection}, $Q_\E$ projects to $d_\M$, then it restricts to a vector field of $\E^*|_{\L_L}$
and gives rise to a connection on $\E|_{\L_L}$ over the $NQ$-manifold 
$(L[1],d_L)$ associated to the Lie algebroid $L$. Its curvature is the restricted curvature of $Q_\E$. 

%Since $L[1]\sim\L_L\rightarrow T^*[2]T[1]M\rightarrow (T\oplus T^*)[1]M$, 
Let us discuss now torsion. We have that 
$p^*L[1]$ is a subbundle of $p^*\Tan[1]M|_{\L_L}$ and it is not difficult to realize that
$\taum|_{\L_L}$ coincides with $\tau_{L[1]}$ under this identification. Let us suppose now that $Q_\E$ restricts to a connection $Q_{L[1]}$ on 
$p^*L[1]$ as in Definition \ref{conn_sub_bundle}; as a consequence, $T_{Q_\E}|_{\L_L}=T_{Q_{L[1]}}$. 

On the other hand, the $K$-curvature and $K$-torsion do not automatically restrict to the curvature and torsion of $Q_{L[1]}$ due to the terms proportional to the coordinate $p_\mu$ measuring 
the difference with the full curvature and torsion in  (\ref{Full_vs_Kcurv}) and (\ref{Full_vs_Ktor}).  We are going to discuss the compatibility conditions between $K$ and $L$ that make these terms vanish so that the $K$-tensors coincide with the full curvature and torsion. 
%We are going to discuss the obstruction on the auxiliary connection $K$ for the $K$-torsion and $K$-curvature to restrict to the Lie algebroid torsion and curvature on the Dirac structure.  
Using (\ref{diracone}), we get
\begin{equation}\label{local_phi_K_L}
  \tilde p_\mu|_{\L_L}
  =\, \frac{1}{2}\Bigl(\phi_{\mu A B} + \KKK{\nu}{\kappa}{\mu} \,\rho^\nu _A \rho_{\kappa B} - \KKK{\nu}{\kappa}{\mu} \,\rho^\nu _B \rho_{\kappa A}\Bigr)\lambda^A\lambda^B= \frac{1}{2}\phi^{KL}_{\mu AB}\lambda^A\lambda^B\;.
\end{equation}
The compatibility is then controlled by the tensor $\phi^{KL}\in\Gamma(\Lambda^2 L^*\otimes T^*M)$ above; in the next section we will give an intrinsic definition of $\phi^{KL}$. 

\begin{prop}\label{obstruction}
Let $D = \nabla+V$ be a generalized connection on $E$ and $\nabla^K$ be an affine connection. Let $L$ be a Dirac structure. The $K$-curvature $R^K_D$ defined in 
(\ref{K_curvature}) restricts to the curvature of $D|_L$ if and only if
$$
\langle V,\phi^{KL}\rangle =0
$$
where $\phi^{KL}$ is defined in (\ref{local_phi_K_L}). If $E=\Tan M$ and $D$ preserves $L\subset\Tan M$ then the $K$-torsion $T^K_D$ defined in (\ref{K_torsion}) 
restricts to the torsion of $D|_L$ if and only if $\phi^{KL}=0$.
\end{prop}

\begin{rem}{\rm
It is interesting to notice that the tensor $\phi^{KL}$ controls the conditions
%We remark that thanks to the characterization given in Remark (\ref{second_def_K_tensor}) an equivalent characterization of this obstruction can be obtained by looking at conditions on $K$ 
that allow
the contraction vector field $\iota^K_a$ defined in (\ref{contraction_vf}) to restrict to a vector field on $\L_L$ for each $a\in\Gamma(L)$. Due to maximality of $\L_L$, it is sufficient to impose the following
\begin{equation}\label{iotacond}
  \omega\Bigl(\phi_* \,\frac{\partial}{\partial \lambda^A}, \iota^K_a\Bigr) = 0\;, \quad \omega\Bigl(\phi_* \frac{\partial}{\partial y^\mu}, \iota^K_a\Bigr) = 0\;,
\end{equation}
where $\omega$ denotes the symplectic form on $T^*[2]T[1]M$ and $\phi_*$ denotes the push forward by the embedding map $\phi : \L_L \hookrightarrow {\cal M}$ as in Appendix \ref{appB}. The first condition in \eqref{iotacond} reads 
$$
  0 = \,\rho^\mu_A a_\mu + \rho_{\mu A} a^\mu \,, \quad% a = a^\mu b_\mu + a_\mu \psi^\mu \in C^1({\cal M})\;,
$$
and is implied by isotropy of $L$. Using this in the second condition of \eqref{iotacond}, we arrive at
$$
%\begin{equation}\label{conditionphi}
 \phi_{\mu A B} + \rho^\nu_A \KKK{\nu}{\kappa}{\mu}\, \rho_{\kappa B} - \rho_{\nu A} \KKK{\kappa}{\nu}{\mu}\, \rho^{\kappa}_B=0\;,
%\end{equation}
$$
which means $ \phi^{KL}=0$ by (\ref{local_phi_K_L}).}
\end{rem}

\smallskip
Finally we discuss two examples of Dirac structures, the tangent bundle and the graph of a Poisson tensor:

\smallskip
\begin{example}{\rm
Let $L=TM$ so that $\L_L$ is defined by $p_\mu=b_\mu=0$ and is isomorphic to $(T[1]M,d_{dR})$. It is then clear that $\phi^L_K=0$ for every $K$.}
\end{example}
\begin{example}{\rm
Let $L={\rm graph}\pi$, where $\pi$ is a Poisson tensor on $M$; $\L_L$
is defined by $\psi^\mu=\pi^{\mu\nu}b_\nu$ and $p_\mu = -\frac{1}{2}\partial_\mu \pi^{\nu\sigma}b_\nu b_\sigma$ and is isomorphic to $(T^*[1]M,d_\pi)$. 
Here $d_\pi$ is the Poisson-Liechnerowicz differential. Indeed, equation \eqref{diracone} of Appendix \ref{appB} in this case reads 
$\phi_\mu{}^{\nu \rho} = -\partial_\mu \pi^{\nu\rho}$.  Finally the  vanishing of $\phi^{KL}$
gives $\nabla^K \pi = 0$, i.e. the Poisson tensor is preserved 
by the connection $\nabla^K$. %One verifies that $Q_\E$ induces a contravariant connection of the Lie algebroid $T^*_\pi$ associated to 
%$\pi$, and a tedious but straightforward calculation shows that the curvature \eqref{allthecomponents}, restricted to 
%$\psi^\mu = \pi^{\mu \nu}b_\nu$ and $p_\mu = -\frac{1}{2} \partial_\mu \pi^{\nu \sigma} b_\nu b_\sigma$ gives the Lie algebroid 
%curvature of $T^*_\pi$, evaluated in the basis determined by the graph of $\pi$.
}
\end{example}

%By construction $Q_{\E}$ preserves $C(\E)$ the space of functions on $\E$ that are linear in the fibre coordinates. We then define the space of sections 
%$\Gamma(\E)=C(\E)_{lin}$. This must be seen as a replacement of $\Gamma(E\otimes\Lambda T^*M)$ in the ordinary case. 

\section{Comparison with naive torsion and curvature}
\label{ClassicalGeometry}

In the previous section we introduced the $K$-curvature in (\ref{K_curvature}) and the $K$-torsion in (\ref{K_torsion}) of a generalized connection $D$ of the exact Courant 
algebroid. From that derivation it must be clear that the role of the auxiliary connection $\nabla^K$ entering the definitions of the tensors is completely arbitrary. 
Moreover, the $K$-curvature (or torsion) is not the full curvature (or torsion) but the missing terms can be ignored only for certain purposes and under strict conditions, 
as the discussion of Dirac structures in the previous section shows. Nevertheless, these $K$-tensors are useful objects and can be used for comparison with other 
structures that have been studied in the literature. This is what we are going to discuss in this section.

We recall that for Lie algebroids, e.g. Dirac structures in Courant algebroids, torsion and curvature are well defined objects and have the same definition as in 
standard Riemannian geometry. However, for Courant algebroids these definitions do not give tensors any more due to the properties of the Courant bracket. 
We call them \emph{naive} to emphasize this, i.e. the naive curvature is the curvature operator $R_D$ defined in \eqref{gualtieri_curvature_operator} and we define the naive torsion operator by
\begin{equation}\label{naivetorsion}
  T_D(a,b):= \,D_ab - D_ba - [[a,b]]\;, \quad a,b \in \Gamma(\Tan M)\;.
\end{equation}

In the case of torsion, the tensor $\mathfrak{T}_D$  defined in \eqref{gualtieri_torsion} is widely used in the literature.

The following proposition compares the naive curvature \eqref{gualtieri_curvature_operator} with the tensorial $K$-curvature $R^K_D$ defined in (\ref{K_curvature}). 
We recall that $K$ is an affine connection entering the definition of $R^K_D$ and $E$ a vector bundle over $M$.

\begin{prop}\label{Kcurvature_vs_naive}
Let $D=\nabla+V$ be a Courant algebroid connection on $E$. The curvature $R^K_D$ reads
\begin{equation}
R^K_D(a,b)= R_D(a,b) + V_{\tilde K(a,b)}\;,
\end{equation}
where
\begin{equation}\label{def_KX} 
\tilde K(X+\xi,Y+\eta) = [[X+\xi, Y+\eta]]_{sk}^H - [X,Y] + \nabla^K_Y(\xi)-\nabla^K_X(\eta)\;.
\end{equation}  
\end{prop}
The proof is a long but straightforward calculation done by expanding $R_D(a,b)$ and comparison with \eqref{K_curvature}. We remark that in case $D$ 
is an ordinary connection, i.e. for vanishing $V$, $R_D^K$ equals the naive $D$-curvature $R_D$ in \eqref{gualtieri_curvature_operator}. In general $R_D$ is not a tensor and $V_{\tilde K(a,b)}$ 
has the right transformation properties that combine with those of $R_D$ in order to ensure the tensoriality of $R^K_D$. We now turn to the torsion. 

\begin{prop}
\label{final expression_torsion}
Let $D=\nabla+V$ be a Courant algebroid connection on 
$\Tan M$. The $K$-torsion is 
then computed as
\begin{equation}
T^K_D(a,b) = T_D(a,b)  + \tilde{K}(a,b) \;.
\end{equation}
where $T_D$ is the operator defined in \eqref{naivetorsion} and $\tilde{K}$ in (\ref{def_KX}).
\end{prop}

Furthermore we observe that the $\tilde{K}$ defined in (\ref{def_KX}) measures the deviation of the $K$-curvature and torsion from the naive definition. 
In particular, given a Dirac structure $L$ that is preserved by $D$, since both the restricted curvature and the torsion coincide with the corresponding 
algebroid curvature and torsion,  the $K$-curvature and $K$-torsion coincide with the naive ones if and only if 
the obstructions discussed in Proposition \ref{obstruction} vanish. Indeed, one directly checks that $\tilde K(a,b)|_L$, written in components along the basis $ e_A = \rho_A^\mu \partial_\mu+ \rho_{A\mu}d x^\mu$ leads to
\begin{align}\label{coordinatestep}
  \tilde K (e_A, e_B)= \Bigl(\frac{1}{2}\Bigl(&\rho^\mu _A \partial_\nu \rho_{T^*\,\mu B} + \partial_\nu(\rho^\mu_B)\rho_{T^*\,\mu A} - A \leftrightarrow B\Bigr) \nonumber \\
  &+\rho^\mu _A \,\KKK{\mu}{\kappa}{\nu}\,\rho_{T^*\,\kappa B} - \rho^\mu_B\,\KKK{\mu}{\kappa}{\nu}\,\rho_{T^*\,\kappa A}\Bigr)dx^\nu\;,
\end{align}
which is precisely \eqref{local_phi_K_L}. We then have that, the tensor $\phi_K^L\in\Gamma(\Lambda^2 L^*\otimes T^*M)$ whose local expression is given in \eqref{local_phi_K_L} reads for each $a,b\in\Gamma(L)$
\begin{equation}\label{global_def_obstr}
\phi^{KL}(a,b) = \tilde{K}(a,b)=\rho_{T^*}([a,b]) + \nabla^K_{\rho(b)}(\rho_{T^*}(a))- \nabla_{\rho(a)}^K(\rho_{T^*}(b))\;.
\end{equation}

%{\color{red}\bf ATTN: the following is not exact. Condition (\ref{intrinsic_condition}) is equivalent to $\rho^*\tilde{K}%(a,b)=0$, so it is a weaker condition. I wonder if 
%it makes sense to write it.}
%%To conclude this section, we want to give a coordinate-free meaning to \eqref{coordinatestep}. For this we observe %from the definition of $\pi_L$ that $\iota_{\rho(e_A)} \rho_{T^*}(e_B) = \pi_L(e_A,e_B)$. Using this and the form of %$e_A \in \Gamma(L) \subset \Gamma((T+T^*)M)$, $e_A = \rho(e_A) + \rho_{T^*}(e_A)$, we get for the one-form part of the Courant bracket of $e_A, e_B$, evaluated on a vector field $\rho(s), s\in \Gamma(L)$,
%\begin{align}
%  \Bigl([[\rho(e_A)+\rho_{T^*}(e_A),&\, \rho(e_B) + \rho_{T^*}(e_B)]]_{sk} \nonumber \\
%  &- [\rho(e_A),\rho(e_B)]\Bigr)(\rho(s)) = \pi_L(s,[e_A,e_B])\;.
%\end{align}
%Finally, we note that this difference is the part of \eqref{def_KX} free of $K$, hence $\tilde K(e_A, e_B) = 0$ is equivalent to 
%\begin{equation}\label{intrinsic_condition}
%    \Bigl( K_{\rho (e_A)}(\rho_{T^*}(e_B)) -  K_{\rho (e_B)}(\rho_{T^*}(e_A))\Bigr)(\rho(s)) = \,\pi_L(s,[e_A,e_B]_L)\;,
%\end{equation}
%where again, $\pi_L = \rho_{T^*}^* \circ \rho \in \Gamma(\wedge^2 L^*)$ is the two-form characterizing the Dirac structure. This is the interpretation of \eqref{conditionphi}: Antisymmetrizing the auxiliary connection with respecto to its arguments is -- on a Dirac structure -- determined by the $L$-two form $\pi_L$ and the bracket on $L$. 

\section{Ricci tensor and scalar curvature}
\label{Ricci}

In order to extract the tensorial objects out of the curvature $R_{Q_\E}$ and torsion $T_{Q_\E}$  of a connection $Q_\E$ we introduced an arbitrary affine connection $K$. To control this arbitrariness, there are two possible meanigful approaches: considering $K$-invariant geometrical quantities or identifying conditions that lead to a canonical choice of $K$. We have seen for instance the conditions on $K$ under which $K$-curvature and $K$-torsion properly restrict to Dirac structures. 

In this  section  we study the dependence of the Ricci tensor and scalar curvature on the choice of the auxiliary connection $K$. 
This also allows a useful comparison with the physics literature. We will compute the Ricci curvature as the appropriate trace of the $K$-curvature and the scalar curvature as contraction of the Ricci tensor with a generalized metric. For alternative approaches to Ricci and scalar curvature we refer to \cite{Garcia-Fernandez:2013gja, Severa:2016lwc}

Let us recall the standard definition of Ricci tensor first. For a Riemannian manifold $(M,g)$ with metric $g\in \Gamma(S^2T^*M)$ and a trivialization given by local sections $e_\mu \in \Gamma(TM)$ with dual $e^\mu$, the components of the Ricci tensor $\textrm{Ric} = \textrm{Ric}_{\mu\nu}e^\mu\otimes e^\nu$ and scalar curvature $\textrm{Scal} \in C^\infty(M)$ are given by
\begin{equation}
  \textrm Ric_{\mu\nu} = \langle e^\sigma, R(e_\sigma,e_\mu)(e_\nu)\rangle\,,\quad \textrm{Scal} = g^{\mu\nu}\textrm{Ric}_{\mu\nu}\;,
\end{equation}
where $R(e_\mu,e_\nu)$ denotes the standard Riemann curvature and $\langle \cdot, \cdot\rangle$ is the evaluation of forms on vectors. We repeat this construction for the $K$-curvature $R^K_D$ of a generalized connection $D=\nabla+V$ on $\Tan M$. Denote by $E_\alpha$ the local basis $(e_\mu,e^\mu)$ of sections of $\Tan M$, and let $E^\alpha$ denote the dual basis with 
respect to the canonical pairing $\langle \cdot,\cdot \rangle$ in $\Tan M$. We define $\textrm{Ric}^K\in\Gamma(\Tan M\otimes\Tan M)$ as $\textrm{Ric}^K_{\alpha\beta} = \langle E^\gamma, R^K_D(E_\gamma, E_\alpha)(E_\beta)\rangle$. 
%In components {\color{red}ATTN: if possible, I would prefer to avoid the graded geometry notation in this section}
%\begin{equation}
%  \textrm{Ric}^K = \textrm{Ric}^K{}_\mu{}^\alpha\, \psi^\mu \frac{\partial}{\partial s^\alpha} + \textrm{Ric}^K{}^{\mu\alpha}\,b_\mu \frac{\partial}{\partial s^\alpha}\;.
%\end{equation}
Finally, contracting with a generalized metric $G \in \Gamma(S^2(\Tan M))$ of the form \eqref{GeneralizedMetric} gives the scalar curvature $\textrm{Scal}^K= G^{\alpha\beta}\textrm{Ric}^K_{\alpha\beta}$. Observing that $E^\alpha$ is identified by $(e^\mu,e_\mu)$, i.e. the index ${}^\alpha$ has contributions ${}^\nu$ and ${}_\nu$, the result consists of four parts
\begin{equation}\label{scalK}
  \textrm{Scal}^K = \textrm{Ric}^{K}{}^{\mu\nu} G_{\mu\nu} + \textrm{Ric}^K{}^\mu_\nu G^\nu _\mu + \textrm{Ric}^K{}^\nu_\mu G^\mu _\nu + \textrm{Ric}^K{}_{\mu\nu}G^{\mu\nu}\;,
\end{equation}
where $G$ has the index structures corresponding to the block matrix form (\ref{GeneralizedMetric}). For convenience of the reader we write here the local coordinate expression of the 
$K$-curvature that can be obtained from  \eqref{allthecomponents}, \eqref{F_explicit_non_graded_vb11_bis} together with (\ref{Full_vs_Kcurv})
\begin{align}
  \cur_{K}{}^\alpha{}_\beta =\, & (\partial_\mu \Gamm{\nu}{\alpha}{\beta}+
\Gamm{\mu}{\alpha}{\gamma}\Gamm{\nu}{\gamma}{\beta} + \tfrac{1}{2}H_{\rho \mu\nu}\,V^\rho{}^\alpha{}_\beta)\psi^\mu\psi^\nu -\VVV{\mu}{\alpha}{\gamma}\VVV{\nu}{\gamma}{\beta}\, b_\mu b_\nu \nonumber \\
  & +(\partial_\mu \VVV{\nu}{\alpha}{\beta}+ \VVV{\nu}{\gamma}{\beta} \Gamm{\mu}{\alpha}{\gamma} -\Gamm{\mu}{\gamma}{\beta} 
\VVV{\nu}{\alpha}{\gamma}-V^{\kappa}{}^\alpha{}_\beta\, K_\mu{}^\nu{}_\kappa)\psi^\mu b_\nu ~.
\nonumber 
  \end{align}
In order to compare the resulting expression for the scalar curvature with the literature, we have to make assumptions on the connection components $V^\mu {}^\alpha{}_\beta$ and $\Gamma_\mu{}^\alpha{}_\beta$. We use the notation introduced in Appendix \ref{appA}.

We impose that the $K$-torsion of $D$ vanishes. 
%The $K$-torsion tensor $T_D^K$ contains the standard torsion as well as additional parts coming from the $V$- as well as $K$-contributions. Setting $T_D^K$ to zero fixes the auxiliary connection in terms of $\Gamma$ and $V$. Taking into account also compatibility with the Courant algebroid pairing $\langle \cdot, \cdot \rangle$,
This fixes the connection $K$ in the following way:
\begin{displaymath}
  K_\mu{}^\nu{}_\rho = \tilde \Gamma_\mu{}_\rho{}^\nu - V^\nu{}_\rho{}_\mu\;,
\end{displaymath}
where $\tilde{\Gamma}_\mu{}_\rho{}^\nu$ denotes the coefficients of the connection $\nabla^{T^*T^*}$ on $T^*M$ as in Appendix \ref{appA}. We further 
impose the connection $D$ to be compatible with the canonical pairing according to (\ref{compatibility_pairing}).
As a consequence, $\nabla^{TT}$ is a torsion-free connection on $TM$ whose dual $(\nabla^{TT})^*$ is identified with $\nabla^{T^*T^*}$. The connection component  $\nabla^{T^*T}$ is seen to be an element of $\Gamma(T^*M\otimes \wedge^2 T^*M$), whose totally antisymmetric part is identified with the three form $H$. We call $\gamma$ the remaining undetermined tensor, i.e. $\gamma$ has components $\Gamma_{\mu\nu\rho}$, symmetric in $(\mu,\rho)$ and antisymmetric in $(\nu,\rho)$. Furthermore, $V^{TT^*}$ is shown to be in $\Gamma(TM\otimes \wedge^2 TM)$ and symmetric in the first and third entry.   
It is easy to see that all other connection components vanish with the only exception of $V^{T^*T} \in \Gamma(TM \otimes \wedge^2 T^*M)$. Hence, the connection $D$ takes the simple form:
\begin{equation}\label{torsion_less_connection}
  D = \nabla + V = \begin{pmatrix}
    \nabla^{TT} & 0\\
    \tfrac{1}{2}\,H +\gamma  & (\nabla^{TT})^*
  \end{pmatrix} + \begin{pmatrix}
    0 & V^{TT^*} \\
    V^{T^*T} & 0
  \end{pmatrix}\,.
\end{equation}
With \eqref{torsion_less_connection}, the scalar curvature takes the form
\begin{align}\label{endscalar}
  \textrm{Scal}(D) =\;& \textrm{Scal}(g,\nabla^{TT}) + G_{\mu\nu}(\nabla_\rho V^{\mu\rho\nu} + V^\sigma{}_\sigma{}_\rho V^{\mu\rho \nu})\nonumber \\
  &- G^{\mu\nu}(\nabla_\mu V^\sigma{}_{\sigma}{}_\nu + V^\sigma{}_\rho{}_\mu V^\rho{}_\sigma{}_\nu)\;.
\end{align}
We denoted by $\textrm{Scal}(g,\nabla^{TT})$ the combination
\begin{displaymath}
  \textrm{Scal}(g,\nabla^{TT}) = G^{\mu\nu}(\partial_\rho \Gamma_\mu{}^\rho{}_\nu - \partial_\mu \Gamma_\rho{}^\rho{}_\nu + \Gamma_\rho{}^\rho{}_\sigma \Gamma_\mu{}^\sigma{}_\nu - \Gamma_\mu{}^\rho{}_\sigma \Gamma_\rho{}^\sigma{}_\nu)\;,
\end{displaymath}
and we notice $G^{\mu\nu} = g^{\mu\nu}$ from \eqref{GeneralizedMetric}. The expression \eqref{endscalar} can be used as a starting point to compare with the literature e.g. in string theory. We remark that to further constrain $V$, metricity of the connection might play a role. This needs a full understanding of the notion of generalized metric in graded geometry which we leave for a future publication.

\subsection*{Acknowledgements}
We are grateful to Christian S{\"a}mann for many discussions at an early stage of this research. P.A. and A.D. thank the university of Florence for hospitality. F.B. and A.D. thank the university of Torino for hospitality. The research of A.D. was supported by OP RDE project No.CZ.02.2.69/0.0/0.0/16\_027/0008495, International Mobility of Researchers at Charles University, as well as by the COST action MP 1405 Quantum structure of spacetime and the Corfu Summer Institute 2019 at EISA.

\bigskip
\noindent {\bf On behalf of all authors, the corresponding author states that there is no conflict of interest}.

\appendix

\section{Generalized connection on $\Tan[1]M$ and torsion components}
\label{appA}
In this appendix we collect our conventions on the component form of a generalized connection on $E = \Tan[1]M$ which are used in the main text, in particular in section \ref{GradedCourant} for the torsion. A generalized connection $D$ is split into two parts, one along a vector field and another along a one form. More precisely, $D = \nabla + V$ where
\begin{displaymath}
\nabla = \left( \begin{array}{cc} \nabla^{TT} & \nabla^{TT^*}\cr
\nabla^{T^*T}&\nabla^{T^*T^*}
\end{array}\right)\;,
\end{displaymath}
where $\nabla^{TT}$ and $\nabla^{T^*T^*}$ are the connections on $TM$ and $T^*M$ defined by the coefficients\footnote{For a convenient track of the type of indices we temporarily use e.g. $\Gamma_\mu{}^{s_\alpha}{}_{s_\beta}$ to denote the component contracting with $s_\alpha \frac{\partial}{\partial s_\beta}$ in the definition of the covariant derivative, where $\alpha$ can be $\mu$ up and $\mu$ down.} 
$\Gamma_\mu{}^\nu{}_\rho\equiv\Gamma_\mu{}^{s_\nu}{}_{s_\rho}$ and $\tilde{\Gamma}_{\mu}{}_\nu{}^\rho \equiv\Gamma_{\mu}{}^{s^\nu}{}_{s^\rho}$, respectively, while $\nabla^{TT^*}\in\Gamma(T^*M\otimes TM^{\otimes 2})$ and $\nabla^{T^*T}\in\Gamma(T^*M{}^{\otimes 3})$ are
defined by $\Gamma_{\mu}{}^\nu{}^\rho\equiv \Gamma_{\mu}{}^{s_\nu}{}_{s^\rho}$ and $\Gamma_{\mu}{}_{\nu}{}_\rho\equiv\Gamma_{\mu}{}^{s^\nu}{}_{s_\rho}$.
Analogously, the component $V$ decomposes as
\begin{displaymath}
V = \left( \begin{array}{cc} V^{TT} & V^{TT^*}\cr
V^{T^*T}&V^{T^*T^*}
\end{array}\right)\;,
\end{displaymath}
so that we define the components of $V^{T^*T}$ by $V^{\mu}{}_{\nu}{}_\rho \equiv V^{\mu}{}^{s^\nu}{}_{s_\rho}$, respectively of $V^{T^*T^*}$ by $V^{\mu}{}_\nu{}^\rho \equiv V^{\mu}{}^{s^\nu}{}_{s^\rho}$. Similarly, for $V^{TT}$ by 
$\tilde{V}^{\mu}{}^\nu{}_\rho \equiv V^{\mu}{}^{s_\nu}{}_{s_\rho}$ and for $V^{TT^*}$ by  $V^{\mu}{}^\nu{}^{\rho}\equiv V^{\mu}{}^{s_\nu}{}_{s^\rho}$.

With these conventions, the direct evaluation of the torsion $T(Q_{\E}) = Q_{\E}(\taum)$ gives, as mentioned in chapter \ref{GradedCourant}, three contributions
\begin{equation}
\label{torsion_components2}
T(Q_\E) = T(\Gamma) + T(V) + T^{(1,1)}(\Gamma,V)
\end{equation}
which we now give explicitely in components. They are
\begin{displaymath}
T(\Gamma)= \Gamma_\mu{}^\rho{}_\nu \psi^\mu\psi^\nu s_\rho + (\tfrac{1}{2}H_{\rho\mu\nu} + \Gamma_{\mu\rho\nu})\psi^\mu\psi^\nu s^\rho 
\in \Gamma(\Lambda^2 T^*M\otimes \Tan M)\,,
\end{displaymath}
\begin{displaymath}
T(V) = V^{\mu}{}^\rho{}^\nu b_\mu b_\nu s_\rho + V^\mu{}_\rho{}^\nu b_\mu b_\nu s^\rho 
 \in \Gamma(\Lambda^2 TM\otimes \Tan M),
\end{displaymath}
and
\begin{equation}\label{torsion_mixed_term}
 T^{(1,1)}(\Gamma,V)=
(\tilde V^{\mu}{}^\rho{}_\nu - \Gamma_\nu{}^\rho{}^\mu)b_\mu \psi^\nu s_\rho + (V^\mu{}_\rho{}_\nu - \tilde \Gamma_\nu{}_\rho{}^\mu)b_\mu \psi^\nu s^\rho +
 p_\mu s^{\mu}\;. 
\end{equation}

It is clear that the first term appearing in $T(\Gamma)$ is the usual torsion
of the affine connection $\nabla^{TT}$. As in the discussion of the $(1,1)$ component of the curvature, we can introduce an affine 
connection on TM in order to get a covariant $\tilde{p}$ so that $T^{(1,1)}(\Gamma,V)-\tilde{p}_\mu s^\mu \in \Gamma(T^*M\otimes TM\otimes \Tan M)$. This gives the $K$-dependent form of the $(1,1)$ component of the torsion, as described and used in \eqref{tensorial_torsion_3} of the main text.

\section{Lagrangian submanifolds of $T^*[2]T[1]M$}\label{Dirac_structures_lagrangian}
\label{appB}
It is known that lagrangian submanifolds of $T^*[2]T[1]M$ invariant under the homological vector field correspond to Dirac structures of the Courant algebroid 
associated to it (e.g. section 4 of \cite{sometitle}). In this appendix we work out explicitely the details of the characterization of Dirac structures in this language 
which is particularly useful for the main text. Let $L \subset \Tan M$ be a Dirac structure; we are going to describe the corresponding $d_\M$-invariant lagrangian submanifold 
${\cal L}_L \stackrel{\phi}{\hookrightarrow} T^*[2]T[1]M$.

As a Lie algebroid, $L$ is described by the $NQ$-manifold $L[1]$ of degree 1. Locally, let us choose for the latter coodinates 
$(y^\mu,\lambda^A)$ for $L[1]$ of degrees $(0,1)$ and for the Courant algebroid $(x^\mu, \psi^\mu, b_\mu, p_\mu)$ as in the main text. Finally let us denote by 
$\omega = dx^\mu dp_\mu + d\psi^\mu db_\mu$ the canonical symplectic structure of $T^*[2]T[1]M$. We describe the embedding $\phi$ by the 
following formulas:
\begin{equation}\label{Dirac_obstr}
 y^\mu= x^\mu,\;\;\; \psi^\mu = \rho^\mu_A \lambda^A\;,\quad b_\mu = \rho_{\mu A} \lambda^A\;, \quad p_\mu = \frac{1}{2}\phi_{\mu A B}\lambda^A\lambda^B \,,
\end{equation}
where $\rho^\mu _A$ and $\rho_{\mu A}$ denote the components of the maps $\rho_{T^*}$ and $\rho$ introduced in section \ref{mainGraded} and $\phi_{\mu AB}$ to be 
determined. %The embedding $\iota$ describes the coordinate change between $(y^\mu, \lambda^A)$ and $(x^\mu,\psi^\mu, b_\mu, p_\mu)$. 
For a basis of vector fields, using \eqref{Dirac_obstr} this means
\begin{align}
  \phi_*(\frac{\partial}{\partial y^\mu}) =\,&\frac{\partial}{\partial x^\mu} + \frac{\partial \rho^\kappa _A}{\partial x^\mu} \lambda^A \frac{\partial}{\partial \psi^\kappa} + \frac{\partial \rho_{\kappa A}}{\partial x^\mu}\lambda^A \frac{\partial}{\partial b_\kappa} + \frac{\partial \phi_{\kappa A B}}{\partial x^\mu} \lambda^A\lambda^B \frac{\partial}{\partial p_\kappa}\,, \\
   \phi_*(\frac{\partial}{\partial \lambda^A}) =\,& \rho^\mu _A \frac{\partial }{\partial \psi^\mu} + \rho_{\mu A}\frac{\partial}{\partial b_\mu} + \phi_{\mu A B}\lambda^B \frac{\partial }{\partial p_\mu}\;.
    \end{align}
${\cal L}_L$ being lagrangian means $\phi^* \omega = 0$, which gives the following two conditions:
\begin{align}
  0=\,& (\phi^*\omega)\Bigl(\frac{\partial}{\partial y^\mu}, \frac{\partial}{\partial \lambda^A}\Bigr) = \,\Bigl(\phi_{\mu A B} - \rho^\kappa _A \frac{\partial \rho_{\kappa B}}{\partial x^\mu} - \frac{\partial \rho^\kappa _B}{\partial x^\mu} \rho_{\kappa A}\Bigr)\lambda^B\,, \label{diracone}\\
  0=\,& (\phi^*\omega)\Bigl(\frac{\partial}{\partial \lambda^A},\frac{\partial}{\partial \lambda^B}\Bigr) = \,\rho^\mu_A \rho_{\mu B} + \rho^\mu_{B}\rho_{\mu A} = 0\;.\label{diractwo}
\end{align}
The first condition fixes $\phi_{\mu AB}$ in the transformations (\ref{Dirac_obstr}) and the second one is satisfied thanks to the fact that $L$ is lagrangian.

It is now a direct computation the check that $d_\M$ restricts to $\L_L$,  {\it i.e.} $\phi_*(d_{L[1]})=d_\M$ is equivalent to involutivity of $L$.

\end{document}